\documentclass[12pt,twoside]{article}

\usepackage[utf8]{inputenc}
\usepackage{amsmath,amsfonts,amssymb}
\usepackage{authblk}

\usepackage{fine-paper}
\usepackage{booktabs}
\usepackage{graphicx}
\usepackage{siunitx}
\usepackage[labelformat=simple]{subcaption}
\usepackage{caption}

\usepackage{nicefrac}
\usepackage[super]{nth}
\usepackage{tikz}
\usepackage{pgfplots}
\usepackage{amsthm}
\newtheorem*{remark}{Remark}

\newcommand{\bS}{\boldsymbol{S}}			% morphology shape tensor
\newcommand{\bu}{\boldsymbol{u}}			% velocity vector
\newcommand{\bE}{\boldsymbol{E}}			% strain rate tensor
\newcommand{\bW}{\boldsymbol{W}}			% vorticity tensor
\newcommand{\II}{I\! I}					% second invariant
\newcommand{\III}{I\! I\! I}				% third invariant

\newcommand{\diff}{\mathrm d}

\newcommand{\etal}{\textit{et~al.}}

\title{Finite-Element Formulation for Advection-Reaction Equations with Change of Variable and Discontinuity Capturing}
\author{S. Haßler$^*$}
\author{A.M. Ranno}
\author{M. Behr}
\affil{Chair for Computational Analysis of Technical Systems (CATS)\\
    Center for Simulation and Data Science (JARA-CSD)\\
    RWTH Aachen University\\
    52056 Aachen, Germany\\
    \null\\
    $^*$hassler@cats.rwth-aachen.de}
\date{}

\begin{document}

\maketitle

\begin{abstract}
We propose a change of variable approach and discontinuity capturing methods to ensure physical constraints for advection-reaction equations discretized by the finite element method. This change of variable confines the concentration below an upper bound in a very natural way. For the non-negativity constraint, we propose to use a discontinuity capturing method defined on the reference element that is combined with an anisotropic crosswind-dissipation operator. This discontinuity capturing cannot completely eliminate negative values but effectively minimizes their occurrence. The proposed methods are applied to different biophysical models and show a good agreement with experimental results for the FDA benchmark blood pump for a physiological red blood cell pore formation model.

\end{abstract}

\section{Introduction}

Advection-diffusion-reaction equations are very common in the course of simulating the generation and distribution of a physical quantity. They can be used, for example, for the determination of a temperature in a fluid, the distribution of a substance that is transported in a fluid (e.g., a drug in blood flow \cite{Zunino2009a}), or even for the estimation of the residence time of a fluid particle in a domain \cite{EsmailyMoghadam2013a}. It is very common that physical constraints exist for the concentration values, such as a non-negativity constraint, or an upper value that cannot be exceeded. However, the numerical schemes cannot always fulfill these constraints, motivating the development of advanced techniques to tackle this problem.

One possible approach is to find a change of variable that confines the quantity of interest by construction to the physical range. Nam \etal~\cite{Nam2011a} proposed a change of variable of quadratic form to tackle the non-negativity constraint for the advection-reaction equation occuring in hemolysis modeling. Ilinca and Pelletier~\cite{Ilinca98a} proposed a logarithmic transformation to ensure the positivity preservation of turbulence variables. And Kramer and Willcox~\cite{Kramer2019a} used change of variables to find a set of equations with a better numerical behavior in the course of model order reduction.

Another possibility is to introduce operators on the partial differential equation (PDE) level or on the discretized system of equations. Hughes and Mallet~\cite{Hughes86e} and Shakib \etal~\cite{Shakib91b} introduced an additional, residual-based, isotropic discontinuity capturing (DC) operator on the PDE level to reduce spurious oscillations at steep internal layers caused by the SUPG method. Codina~\cite{Codina93a} has a similar approach but proposed a crosswind-dissipation (CWD) DC method that only acts perpendicular to the streamlines. As an example of a method that acts on the discretized system of equations, Kuzmin and Turek~\cite{Kuzmin2004a} use a total variation diminishing for the finite element method that adds a minimum amount of diffusion to ensure a positive solution.

Ventricular assist devices (VADs) are mechanical heart pumps that support a failing heart and can be used as a bridge-to-transplant or as a destination therapy \cite{Kirklin2017a}. During the development of such devices, computational fluid dynamics (CFD) is a valuable assessment tool. While the prediction of flow characteristics and hydraulic performance is now sufficiently reliable, the hemolysis estimation remains a challenging task. Hemolysis, the release of hemoglobin from the red blood cells (RBCs) to the blood plasma, can be modeled with different complexities. One of the simplest approaches is the so-called stress-based model, which determines the released hemoglobin depending on a scalar shear stress and an exposure time to that stress, using a power law fitted to experimental data \cite{Giersiepen90a,Zhang2011a,Ding2015a}.

Further developments to more complex models led to the so-called strain-based approach that relies on an intermediate model to predict the cell deformation of the RBCs. Arora \etal~\cite{Arora2004a} proposed the morphology model, a droplet-like model to compute the ellipsoidal elongation and rotation of RBCs in blood flow. The model computes an effective shear stress acting on the RBC, which is then used in the power law to predict the hemolysis caused by the flow. Vitale \etal~\cite{Vitale2014} derived a pore formation model that computes the hemoglobin release through pores forming on the stretched surface, which is computed with the morphology model. The threshold model of Chen and Sharp~\cite{Chen2011a} and the viscoelastic model of Arwatz and Smits~\cite{Arwatz2013a} are further examples for strain-based models. New developments towards the correct description of the RBC deformation were proposed by Ezzeldin \etal~\cite{Ezzeldin2015} and Sohrabi and Liu~\cite{Sohrabi2017a}. For more details on hemolysis modeling, we refer to two recent review articles by Yu \etal~\cite{Yu2017a} and Faghih and Sharp~\cite{Faghih2019a}.

The next section will start with an introduction of three different hemolysis modeling approaches: stress-based and strain-based models dependent of the empirical power law, and the physiologically motivated RBC pore formation model. In section 3, we discuss a change of variable approach and discontinuity capturing methods for a general advection-reaction equation to suppress nonphysical values. The proposed techniques are tested with three different test cases in section 4: a two-dimensional channel test case, an academic two-dimensional blood pump test case, and the benchmark VAD proposed by the U.S. Food and Drug Administration.

\section{Advection-Reaction Equations}

In this article, we will use pure advection-reaction equations for the quantification of flow-induced mechanical hemolysis. Hemolysis is the release of RBC enclosed hemoglobin to the blood plasma through pores that can form on a stressed red blood cell membrane, or due to the complete rupture of highly stretched RBCs. In the following subsections, we will introduce the empirical power law model and the physiologically inspired pore formation model for the estimation of the release and distribution of hemoglobin in the blood plasma.

\subsection{Power Law Model for Hemolysis Estimation}

The power law model is based on empirical findings dating back to Blackshear and Blackshear~\cite{Blackshear87a} and Leverett \etal~\cite{Leverett72a} who found that the mechanical hemolysis strongly depends on the fluid shear stresses $\sigma_s$ and the exposure time $\tau$ to these stresses. Experiments in Couette shearing devices conducted by Heuser and Opitz~\cite{Heuser80a} and Wurzinger \etal~\cite{Wurzinger86a} supported this relation and Giersiepen \etal~\cite{Giersiepen90a} proposed a power law
\begin{equation}
IH = A \sigma_s^\alpha \tau^\beta
\end{equation}
to compute the index of hemolysis $IH$, the ratio of plasma-free hemoglobin to total hemoglobin content, based on Wurzinger's experimental data.

Song \etal~\cite{Song2003a} derived a different parameter set based on the experimental data of Heuser \etal~\cite{Heuser80a}. However, Paul \etal~\cite{Paul2003a} concluded that these early experiments contained a high amount of secondary hemolysis caused by the experimental devices. Recent hemolysis experiments were conducted by Zhang \etal~\cite{Zhang2011a} and Ding \etal~\cite{Ding2015a} for different blood species and corresponding power law parameters were proposed. An overview of the parameter sets is shown in Table~\ref{tbl:Power-Law-Parameters}.
\begin{table}[b]
\begin{center}
\caption{Power-law parameters derived by different research groups.}
\begin{tabular}{l c c c c}
\toprule
Research group                        &       $A$        & $\alpha$ & $\beta$ & blood species \\
\midrule
Giersiepen \etal~\cite{Giersiepen90a} &  $\num{3.62e-7}$ &  2.416   &  0.785  &     human     \\
Song \etal~\cite{Song2003a}           &   $\num{1.8e-8}$ &  1.991   &  0.765  &    porcine    \\
Zhang \etal~\cite{Zhang2011a}         & $\num{1.228e-7}$ &  1.9918  &  0.6606 &     ovine     \\
Ding \etal~\cite{Ding2015a}           & $\num{3.458e-8}$ &  2.0639  &  0.2777 &     human     \\
Ding \etal~\cite{Ding2015a}           & $\num{6.701e-6}$ &  1.0981  &  0.2778 &    porcine    \\
\bottomrule
\end{tabular}
\label{tbl:Power-Law-Parameters}
\end{center}
\end{table}

The numerical integration of the power law in space and time for complex flows can be done in an Eulerian or Lagrangian approach. For an overview on Lagrangian approaches to hemolysis modeling we refer to Yu \etal~\cite{Yu2017a}. We will focus on the Eulerian formulation which results in an advection-reaction equation. Garon and Farinas~\cite{Garon2004a} and Farinas \etal~\cite{Farinas2006a} proposed to use the power law linearized in (exposure) time to determine the reaction term, yielding
\begin{equation}
\left(\frac{\partial}{\partial t} + \bu \cdot \nabla \right) l_{IH} = \left(A \sigma_s^\alpha\right)^{1/\beta} \left(1 - l_{IH}\right),
\label{eq:AR-power-law}
\end{equation}
with the linearized index of hemolysis $l_{IH} = IH^{1/\beta}$ and the flow velocity $\bu$.

\enlargethispage{6mm}What is left to complete the model formulation is to define the scalar shear stress $\sigma_s$. We will present two modeling strategies, the stress-based and the strain-based approach, hereafter.

\subsubsection{Stress-Based Hemolysis Model}

A first approximation is to use a scalar measure of the flow velocity gradient for $\sigma_s$, assuming an instantaneous deformation of the RBCs to the flow, which is commonly named stress-based model.
One possible choice is to use the second invariant $\II_{\bE} $ of the strain rate tensor $\bE = \left(\nabla \bu + \nabla \bu^T \right) / 2$ for the computation of the shear stress
\begin{equation}
\sigma_s = 2 \mu \sqrt{-\II_{\bE}},
\end{equation}
with the blood viscosity $\mu$, which is then inserted in eq.~\eqref{eq:AR-power-law}.

\subsubsection{Strain-Based Hemolysis Model}

Arora \etal~\cite{Arora2004a} proposed a droplet-like model to first compute a time-dependent, effective shear stress acting on the RBCs. This effective shear stress can then be used as a scalar measure for $\sigma_s$, which is called a strain-based approach.

The droplet-like model, called morphology model, describes the RBCs by an ellipsoidal shape tensor $\bS$. While Arora \etal~\cite{Arora2004a} used a Lagrangian description, Pauli \etal~\cite{Pauli2013a} reformulated the problem in an Eulerian way. To account for the deformation, rotation and relaxation of RBCs in blood flow $\bu$, the shape tensor has to fulfill the partial differential equation
\begin{equation}
\left( \frac{\partial}{\partial t} + \bu\cdot\nabla \right) \bS
= \underbrace{-\alpha_{1}\left(\bS-g\!\left(\bS\right)\boldsymbol{1}\right)}_{\textrm{relaxation}}
+ \underbrace{\alpha_{2}\left(\bE\bS+\bS\bE\right)}_{\textrm{elongation}}
+ \underbrace{\alpha_{3}\left(\bW\bS-\bS\bW\right)}_{\textrm{rotation}},
\label{eq:Res-morph}
\end{equation}
with unit matrix $\boldsymbol{1}$, the vorticity tensor $\bW = \left( \nabla\bu - \nabla\bu^T \right) / 2$, problem specific parameters $\alpha_1 = \SI{5}{s^{-1}}$ and $\alpha_2 = \alpha_3 = \num{4.2298e-4}$ (cf. \cite{Arora2004a}), and a scalar $g(\bS) = 3 \III_{\bS} / \II_{\bS}$ to ensure the volume conservation of the RBCs. Recently, we demonstrated how the morphology equation can be numerically stabilized by using a variational multiscale formulation for a logarithmic transformation of the shape tensor $\bS$ \cite{Hassler2019a}.

For a computed shape $\bS$, one can compute the longest and shortest semi-axis of the ellipsoid, $L$ and $W$, with which one can determine the distortion $D = (L - W) / (L + W)$, a measure for the deformation of the RBC. Finally, one can compute an effective shear stress
\begin{equation}
\sigma_s = \sigma_{\mathrm{eff}} = \frac{2 \mu \alpha_1 D}{\left(1 - D^2\right) \alpha_2},
\end{equation}
from the distortion and use it as the scalar shear stress measure in eq.~\eqref{eq:AR-power-law}.

\subsection{Pore Model for Hemolysis Estimation}

Vitale \etal~\cite{Vitale2014} proposed a different approach to hemolysis estimation that takes the physiological processes into account and does not rely on empirical findings as the power law model does. They consider the energies of a stressed RBC membrane and determine the number and size of pores forming in its lipid bi-layer. The RBC deformation needed for this model also comes from the morphology model. The hemoglobin contained inside the RBC can then diffuse through the pores to the blood plasma based on Fick's law. Although Vitale \etal~\cite{Vitale2014} proposed a Lagrangian approach in their article, their model can be easily modified to be used in the Eulerian frame.

As mentioned before, the deformation of the RBC is taken from the solution of the morphology model eq.~\eqref{eq:Res-morph}. But, instead of using the distortion $D$ as a measure for the deformation, the RBC surface area strain
\begin{equation}
\varepsilon = \frac{A_{\bS} - A_0}{A_0}
\end{equation}
is computed, with the surface area of the deformed RBC $A_{\bS}$, and its initial surface area $A_0$. For a given surface strain $\varepsilon$, the number and size\footnote{The model makes the assumption that all pores have the same area.} of the pores, and hence the total pore area $A_p(\varepsilon)$, can be determined by the minimization of the total membrane energy. It is assumed that the membrane energy is mainly driven by two competing factors: the release of stresses when pores are forming and the increase in energy according to the exposure of hydrophobic lipid tails in the RBC lipid bi-layer to the blood plasma. This pore formation process only starts when a critical threshold value $\varepsilon_0 = \SI{0.16}{\%}$ is exceeded. For a detailed discussion we refer to Vitale \etal~\cite{Vitale2014}.

The diffusion of hemoglobin through the pores is driven by a concentration gradient; rearranging Fick's law and using the material derivative yields the advection-reaction equation
\begin{equation}
\left(\frac{\partial}{\partial t} + \bu \cdot \nabla\right) IH = \kappa \frac{A_p(\varepsilon)}{V_{\mathrm{RBC}}} \left( 1 - \frac{IH}{1-Hct} \right),
\label{eq:poreModel}
\end{equation}
with a mass transfer coefficient $\kappa$, the total volume of RBCs in blood $V_{\mathrm{RBC}}$, and the hematocrit $Hct$. The mass transfer coefficient can be modeled, dependent on the fluid shear rate $G_f$, as
\begin{equation}
\kappa = h G_f^k
\end{equation}
with coefficients $(h,k) = (\num{4.48e-8},\num{1.31})$ which were fitted to the porcine blood experiments of Ding \etal~\cite{Ding2015a}.

\section{Circumventing Nonphysical Concentration Values}

The general residual form of an advection-reaction equation with a saturation effect that we yield during the modeling of the generation and transportation of a molecule in blood flow is
\begin{equation}
\mathcal{R}(c) = \left( \frac{\partial}{\partial t} + \bu \cdot \nabla \right) c - \mu \left( \nu - c \right) = 0.
\end{equation}
Herein, $c$ is the concentration of the molecule (e.g., hemoglobin or an antiproliferative drug) in blood, $\bu$ is the blood velocity and $\mu$ and $\nu$ are constants that do not depend on $c$. The corresponding values for different concentration models are shown in Table~\ref{tbl:Conc-Models}.
\begin{table}
\centering
\caption{Different concentration models and their constants.}
\label{tbl:Conc-Models}
\begin{tabular}{l c c c}
\toprule
     Model      &    Molecule     &                        $\mu$                        &   $\nu$   \\
\midrule
Power-Law Model &   Hemoglobin    &           $(A \sigma_s^\alpha)^{1/\beta}$           &     1     \\
Pore Model      &   Hemoglobin    & $\frac{\kappa}{1-Hct}\frac{A_p(\alpha)}{V_{RBC}}$ & $1 - Hct$ \\
Drug Model      & e.g., Everolimus &                          0                          &  $c_S^0$  \\
\bottomrule
\end{tabular}
\end{table}
Naturally, the concentration $c$ is confined between 0 and $\nu$, which might be violated during the numerical simulation, though.

\subsection{Change of Variable for the Concentration}

One possibility to circumvent nonphysical values is to use a change of variable for the concentration. The transformation
\begin{equation}
c = \nu \left( 1 - \exp\!\left(-\frac{\bar{c}}{k}\right) \right), \label{eq:UpperTrans}
\end{equation}
with a scaling constant $k$ restricts $c$ to values smaller than $\nu$ by definition. This yields a transformed advection-reaction equation
\begin{equation}
\mathcal{R}(\bar{c}) = \left(\frac{\partial}{\partial t} + \boldsymbol{u} \cdot \nabla \right) \bar{c} - k \mu = 0.
\label{eq:AR-transformed}
\end{equation}

In principle, one could also think about a transformation of
\begin{equation}
c = \frac{\nu}{1 + \exp\!\left(-\frac{\bar{c}}{k}\right)} \label{eq:FullTrans}
\end{equation}
that confines $c$ between 0 and $\nu$. However, for models where $\mu \neq 0$ this results in a source term in the transformed equation which explicitly contains an exponential $\exp\!\left(-\frac{\bar{c}}{k}\right)$. This leads to numerical problems, since it results in a huge source term for small concentration values. Nevertheless, the transformation might be useful for models without a source term ($\mu = 0$).

\begin{remark}
In the case of a drug model (see Table~\ref{tbl:Conc-Models}), the boundary conditions are defined by the drug release rate into the blood flow. Under certain assumptions \cite{Vergara2008}, the imposition of the release rate at the interface between the stent and the blood flow results in the Robin boundary condition:
\begin{equation*}
D \nabla c \cdot \boldsymbol{n} + \varphi (t) c = \varphi (t) c_S^0,
\end{equation*}
where $D$ is the diffusion parameter, $ \boldsymbol{n}$ the boundary unit outward normal and $\varphi(t)$ is a time dependent function.  This boundary condition ensures the dependency of the released drug $c$ to the initial drug charge of the stent $c_S^0$. Using the change of variable (\ref{eq:UpperTrans}), this dependency is still preserved, although resulting in the Neumann boundary condition:
\begin{equation*}
D \nabla \bar{c} \cdot \boldsymbol{n} =\varphi (t) c_S^0.
\end{equation*}
This is easier to implement and does not require the addition of an unknown surface term integrating the concentration variable $c$ in the following finite element formulation (\ref{eq:WeakAR-transformed}), which the Robin boundary condition entails.
\end{remark}

\subsection{Finite Element Formalism}

Equation \eqref{eq:AR-transformed} is solved using a streamline-upwind/Petrov-Galerkin (SUPG) finite element formulation. The physical domain $\Omega$ is decomposed in a finite element mesh $\overline{\Omega}_n$ with $n$ linear elements $\Omega_n$. To describe a well-formed system, we have to prescribe boundary conditions on the inflow part of the boundary $\Gamma^{\mathrm{in}}$ (which are always Dirichlet boundary conditions in our case) and an initial condition. Let the trial solution and weighting function spaces be denoted by
\begin{align}
\mathcal{S}^h_n = \left\lbrace \bar{c}^h \in C^0\!\left(\overline{\Omega}_n\right) \middle| \left.\bar{c}^h\right|_{\Gamma^\mathrm{in}} = g_c \right\rbrace, \\
\mathcal{V}^h_n = \left\lbrace w^h \in C^0\!\left(\overline{\Omega}_n\right) \middle| \left.w^h\right|_{\Gamma^\mathrm{in}} = 0 \right\rbrace.
\end{align}
The weak form of eq. \eqref{eq:AR-transformed} then becomes: find $\bar{c}^h \in \mathcal{S}_n$ for the given initial condition $\bar{c}^h(t=0) = \bar{c}_0$ such that $\forall w^h \in \mathcal{V}_n$
\begin{equation}
0 = \int_{\Omega_n} w^h \left( \frac{\partial \bar{c}^h}{\partial t} + \bu \cdot \nabla \bar{c}^h - k \mu \right) \diff \Omega 
  + \int_{\Omega_n} \tau \left( \bu \cdot \nabla w^h \right) \mathcal{R}\!\left(\bar{c}^h\right) \diff \Omega
  \label{eq:WeakAR-transformed}
\end{equation}
is satisfied. The stabilization parameter $\tau$ is chosen according to Shakib \etal~\cite{Shakib91b} as
\begin{equation}
\tau = \left( \left(\frac{2}{\Delta t}\right)^2 + \bu \cdot \boldsymbol{G} \bu \right)^{-\frac{1}{2}},
\end{equation}
where $\boldsymbol{G}_{ij} = \sum_k \frac{\partial \xi^k}{\partial x^i} \frac{\partial \xi^k}{\partial x^j}$ is the covariant metric tensor mapping to a symmetric\footnote{We include the mapping to an equilateral triangle or a regular tetrahedron from Pauli and Behr \cite{Pauli2016b} in the definition of the metric tensor.} reference element \cite{Pauli2016b}.

\subsection{Discontinuity Capturing Methods}

If we are in a situation where we cannot make use of a transformation that prohibits negative concentration values, an additional discontinuity capturing (DC) operator is introduced to eq.~\eqref{eq:WeakAR-transformed} to minimize the occurences of negative concentration values. Shakib \etal~\cite{Shakib91b} introduced an isotropic discontinuity capturing operator defined on the reference element as
\begin{equation}
\int_{\Omega_n} \nu_\mathrm{DC}\!\left(\mathcal{R}\!\left(\bar{c}^h\right)\right)\nabla w^h \cdot \boldsymbol{G}^{-1} \nabla \bar{c}^h \diff\Omega,
\end{equation}
and proposed two different definitions for the numerical diffusion $\nu_\mathrm{DC}$, a linear form
\begin{equation}
\nu_\mathrm{DC-lin}\!\left(\mathcal{R}\!\left(\bar{c}^h\right)\right) = \sqrt{\frac{\mathcal{R}\!\left(\bar{c}^h\right)^2}{\nabla \bar{c}^h \cdot \boldsymbol{G}^{-1} \nabla \bar{c}^h}},
\end{equation}
and a quadratic one
\begin{equation}
\nu_\mathrm{DC-quad}\!\left(\mathcal{R}\!\left(\bar{c}^h\right)\right) = 2\frac{\tau \mathcal{R}\!\left(\bar{c}^h\right)^2}{\nabla \bar{c}^h \cdot \boldsymbol{G}^{-1} \nabla \bar{c}^h}.
\end{equation}
The nice thing about this formulation is that the element lengths are intrinsically contained in the covariant and contravariant metric tensors $\boldsymbol{G}$ and $\boldsymbol{G}^{-1}$.

Codina~\cite{Codina93a} proposed to use an anisotropic discontinuity capturing operator that acts only perpendicular to the streamlines, which is called a crosswind-dissipation (CWD) technique. The corresponding term in the weak form is
\begin{equation}
\int_{\Omega_n} \nu_\mathrm{DC}\!\left(\mathcal{R}\!\left(\bar{c}^h\right)\right)\nabla w^h \cdot \left( \boldsymbol{1} - \frac{1}{|\bu|^2} \bu\otimes\bu \right) \nabla \bar{c}^h \diff\Omega.
\end{equation}
His proposed numerical diffusion is given by
\begin{equation}
\nu_\mathrm{Cod} = \frac{1}{2} h^e \max\!\left( 0,C - \frac{1}{Pe^e(\bu_\Vert)} \right) \frac{\left| \mathcal{R}\left(\bar{c}^h\right) \right|}{\left| \nabla \bar{c}^h \right|},
\end{equation}
with the explicit element length $h^e$, an element dependent constant $C$, and the element Péclet number $Pe^e$ of $\bu_\Vert$, the projection of $\bu$ onto $\nabla \bar{c}^h$.

In this article, we propose to use Shakib's DC method in combination with a CWD technique, yielding
\begin{equation}
\int_{\Omega_n} \nu_\mathrm{DC}\!\left(\mathcal{R}\!\left(\bar{c}^h\right)\right)\nabla w^h \cdot \frac{\partial\boldsymbol{x}}{\partial\boldsymbol{\xi}} \left( \boldsymbol{1} - \frac{1}{\bu\cdot\boldsymbol{G}\bu} \bu \otimes \bu \right) \frac{\partial\boldsymbol{x}}{\partial\boldsymbol{\xi}}^T \nabla \bar{c}^h \diff\Omega,
\end{equation}
and using the definitions of $\nu_\mathrm{DC-lin}$ and $\nu_\mathrm{DC-quad}$ for the numerical diffusion. Since we observe convergence problems for the simulations, we further propose a linearisation of the highly non-linear DC operators. To this end, we compute the numerical diffusion $\nu_\mathrm{DC}$ with the concentration values from the previous time step and solve for three steady time steps in total. This technique is applied to the test cases in the following section. Further DC techniques and their corresponding definitions of the numerical diffusion are reviewed in the article of John and Knobloch~\cite{John2007a}.

\section{Results}

In this section, we will apply the previously presented change of variable and discontinuity capturing methods to three different test cases: A two-dimensional, academic power law problem for a flow in a channel, a simple, two-dimensional blood pump for the stress-based hemolysis estimation, and the benchmark blood pump, proposed by the U.S. Food and Drug Administration (FDA), as a validation of the pore model hemolysis prediction. The computations for this section were performed on the supercomputer JURECA at Forschungszentrum Jülich \cite{jureca}.

\subsection{Two-dimensional Channel}

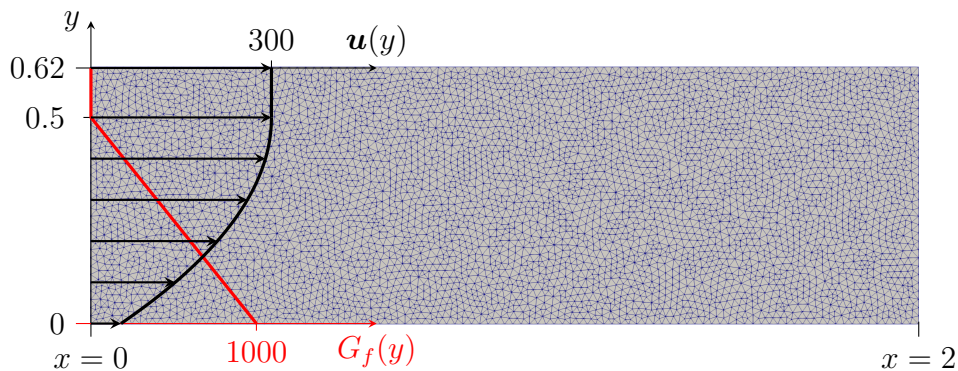
\begin{figure}[b]
\begin{center}
\begin{tikzpicture}[>=stealth]
\node[inner sep=0pt,anchor=south west] at (-0.025,0.715) {\includegraphics[height=3.51cm]{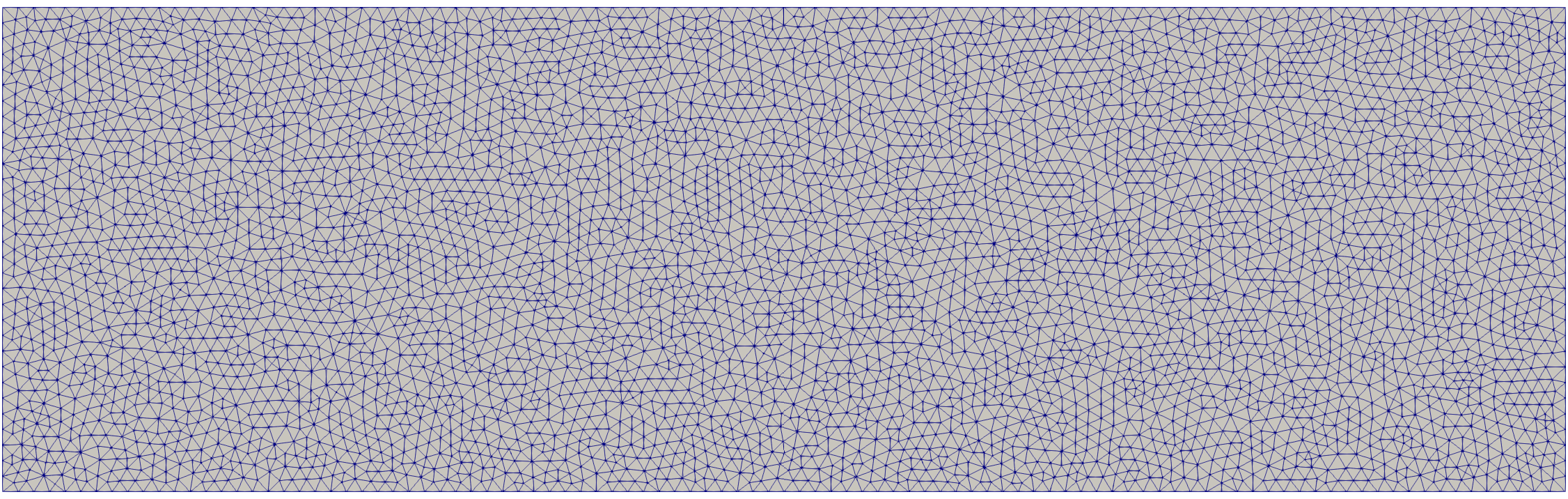}};
\draw[->] (0,0.6) node[below] {$x=0$} -- +(0,4.2) node[left] {$y$};
\draw (11.005,0.8) -- +(0,-0.2) node[below] {$x=2$};
\draw[->] (-0.2,4.175) node[left] {$0.62$} -- +(4,0) node[above] {$\bu(y)$};
\draw[->,red] (-0.2,0.775) -- +(4,0) node[below] {$G_f(y)$};
\node[left] at (-0.2,0.775) {$0$};
\draw (0,3.517) -- ++(-0.1,0) ++(-0.1,0) node[left] {$0.5$};

\draw[very thick,red] (0,4.175) -- ++(0,-0.658) -- ++(2.2,-2.742);
\draw[red] (2.2,0.78) -- +(0,-0.1) node[below] {$1000$};

\begin{scope}[thick,->]
    \draw (0,0.775) -- +(0.4,0);
    \draw (0,1.3234) -- +(1.12,0);
    \draw (0,1.8718) -- +(1.68,0);
    \draw (0,2.4202) -- +(2.08,0);
    \draw (0,2.9685) -- +(2.32,0);
    \draw (0,3.517) -- +(2.4,0);
    \draw (0,4.175) -- +(2.4,0);
\end{scope}
\draw[very thick] (2.4,4.175) -- +(0,-0.658);
\draw[very thick,domain=0.775:3.517,smooth,variable=\y] plot ({2.4-0.266008*(\y-3.517)*(\y-3.517)},{\y});
\draw (2.4,4.175) -- +(0,0.1) node[above] {$300$};
\end{tikzpicture}
\caption{Computational mesh of the two-dimensional channel test case together with the flow profile and the instantaneous shear rate.}
\label{fig:ChannelMeshFlow}
\end{center}
\end{figure}

The first test case that we are investigating is a two-dimensional channel with $n_e = \num{9484}$ unstructured triangular elements and $n_n = \num{4874}$ nodes that is shown in Fig.~\ref{fig:ChannelMeshFlow} (all lengths in \si{cm}).
We consider the power law model with parameters $A = 1.0$, $\alpha = 2.0$, and $\beta = 1.0$,
prescribe a parabolic velocity profile as
\begin{equation}
\bu = u \boldsymbol{e}_x,\quad u = \begin{cases} 300 - 1000 (0.5 - y)^2, & 0 < y < 0.5 \\ 300, & y \geq 0.5 \end{cases},
\end{equation}
and arbitrarily choose a fluid viscosity of $\SI{0.35}{g/cm/s}$. At the inflow, we prescribe a concentration of $c = 0$. This academic problem shows analytical concentration values ranging from 0 to 1 inside the domain.

\begin{figure}
\center\includegraphics[width=.85\textwidth]{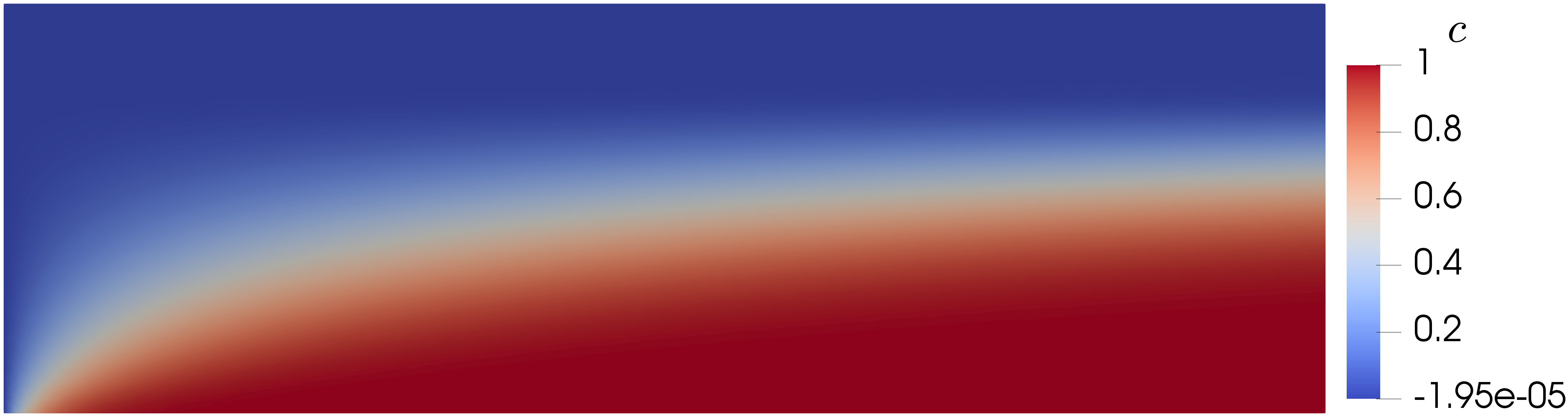}
\caption{Concentration values in the channel test case using the transformation~\eqref{eq:UpperTrans} and a crosswind-dissipation discontinuity capturing method with the quadratic form $\nu_\mathrm{DC-quad}$.}
\label{fig:Conc2dChannel}
\end{figure}

For the untransformed concentration without a DC method we observe nonphysical values of up to $\num{-2.61e-3}$ and $\num{1.002}$. While the transformation $c = 1 - \exp(-\bar{c})$ completely eliminates values above 1, the discontinuity capturing operator with $\nu_\mathrm{DC-quad}$ is able to reduce the negative values by two orders of magnitude (cf. Fig.~\ref{fig:Conc2dChannel}), whereas $\nu_\mathrm{DC-lin}$ even eliminates the negative values completely.

\begin{figure}[b!]
\begin{minipage}{.48\textwidth}
\begin{center}
\begin{tikzpicture}
\node[inner sep=0pt,anchor=south west] at (0,0) {\includegraphics[width=\textwidth]{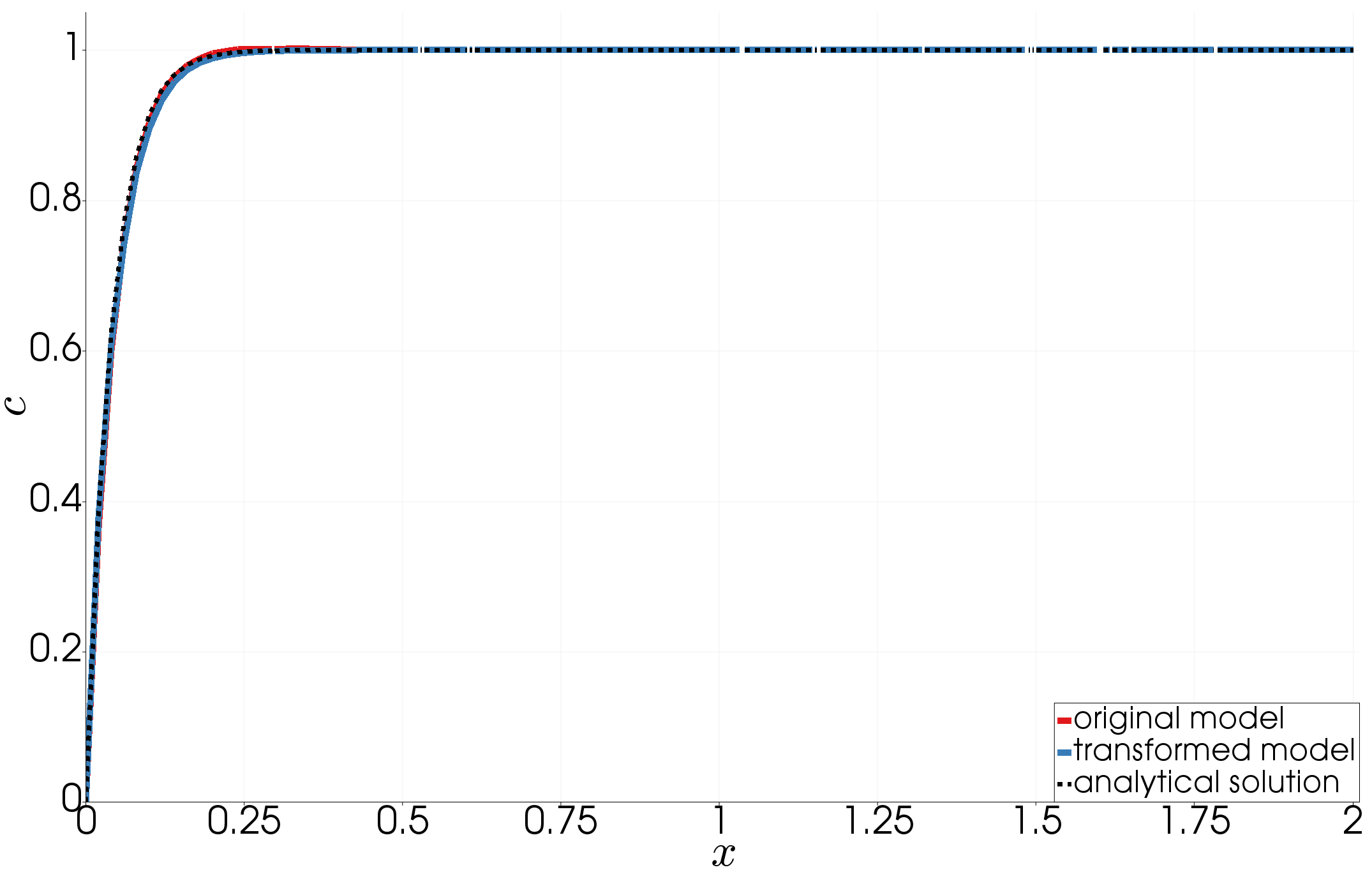}};
\node[inner sep=0pt,anchor=south west] at (1.7,1.1) {\includegraphics[width=.7\textwidth]{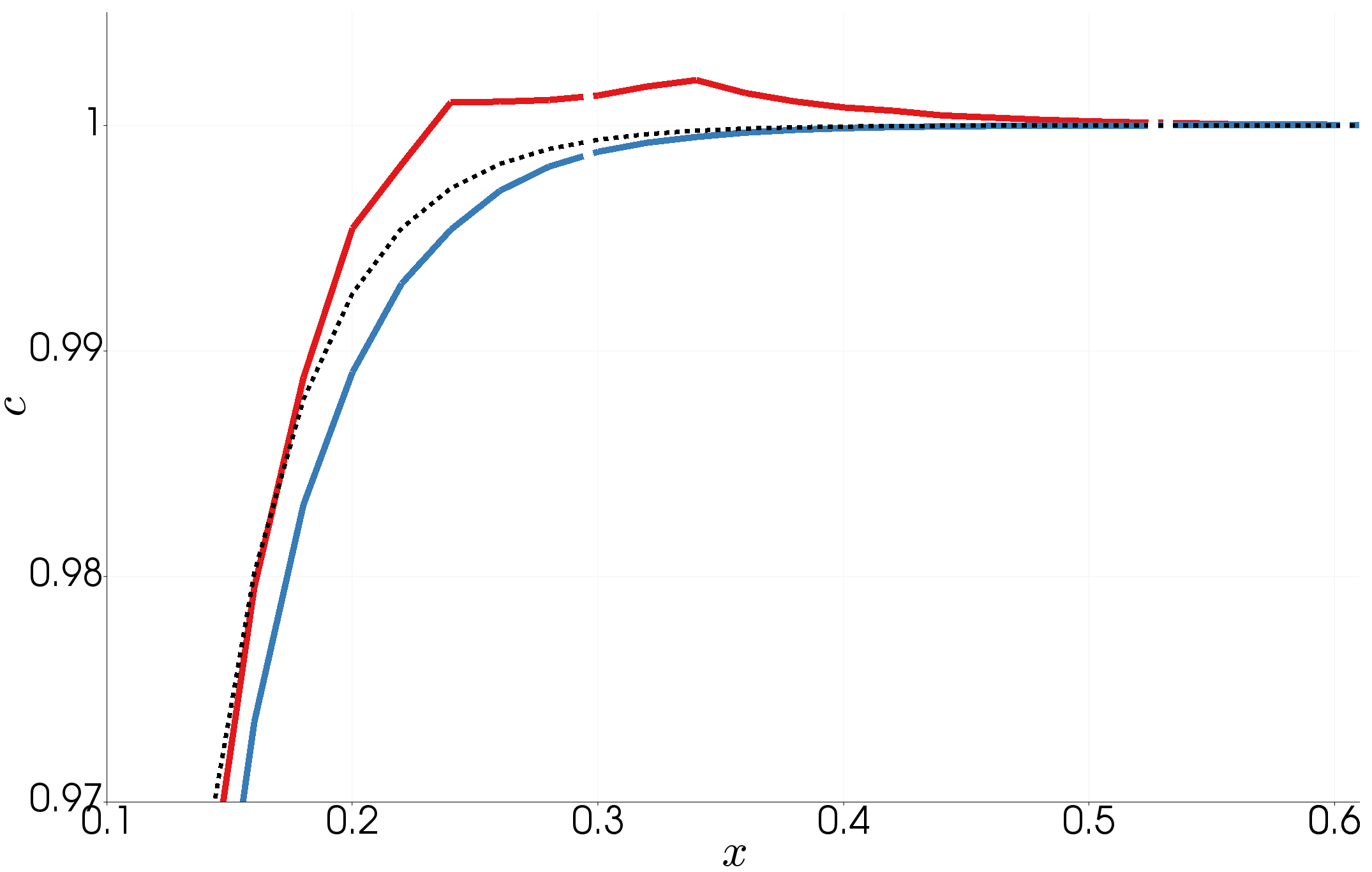}};
\draw[dotted] (0.966,5.1508) rectangle +(1.624,0.25);
\draw[dotted] (0.966,5.1508) -- (2.163,1.49);
\draw[dotted] (2.57,5.4008) -- (7.9,4.97);
\end{tikzpicture}
\captionof{figure}{Comparison of the concentration values on the lower wall for the original and the transformed model.\\\null}
\label{fig:CompWall-Channel}
\end{center}
\end{minipage}\hfill
\begin{minipage}{.48\textwidth}
\begin{center}
\begin{tikzpicture}
\node[inner sep=0pt,anchor=south west] at (0,0) {\includegraphics[width=\textwidth]{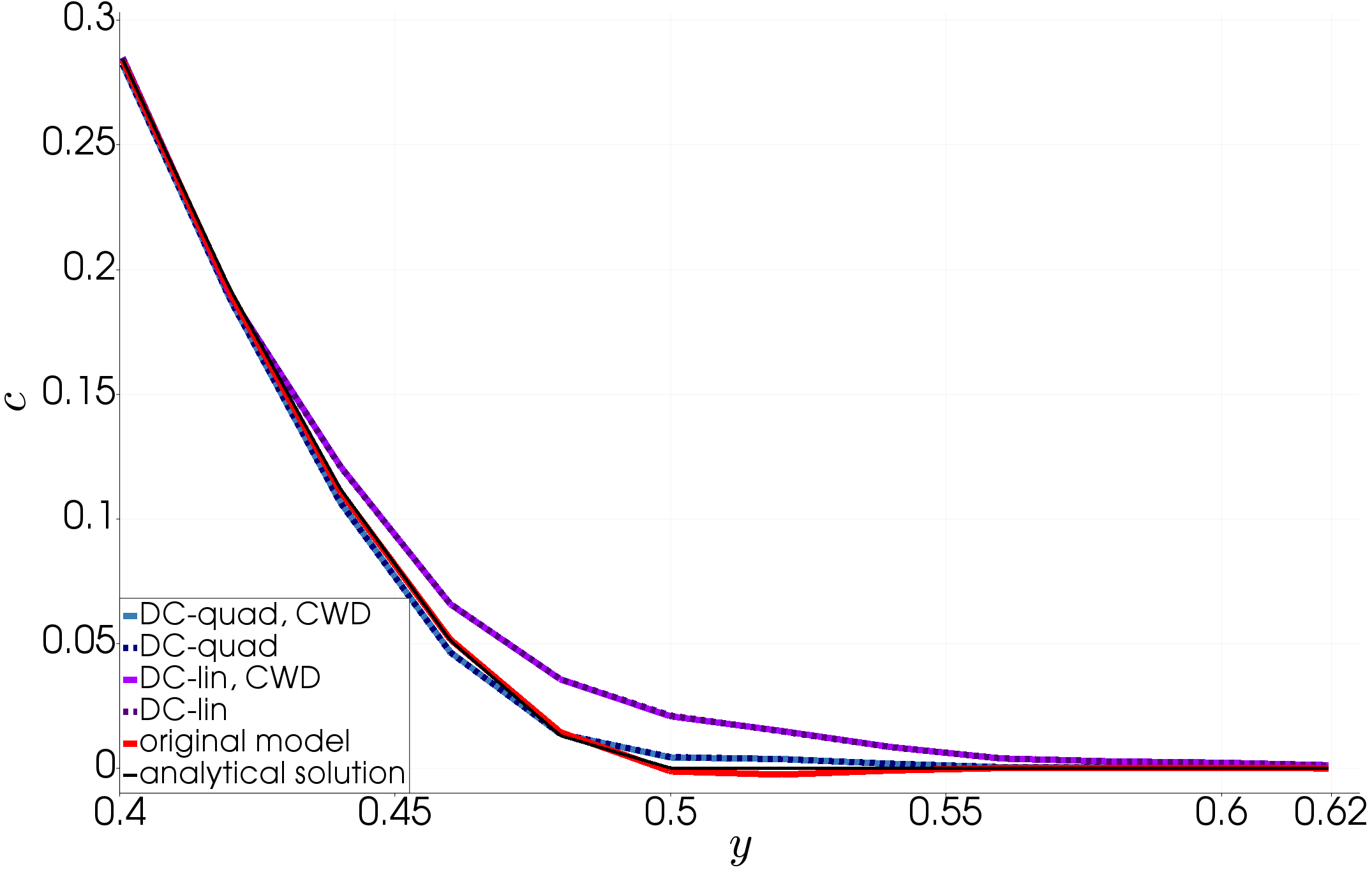}};
\node[inner sep=0pt,anchor=south west] at (2.85,1.85) {\includegraphics[width=.67\textwidth]{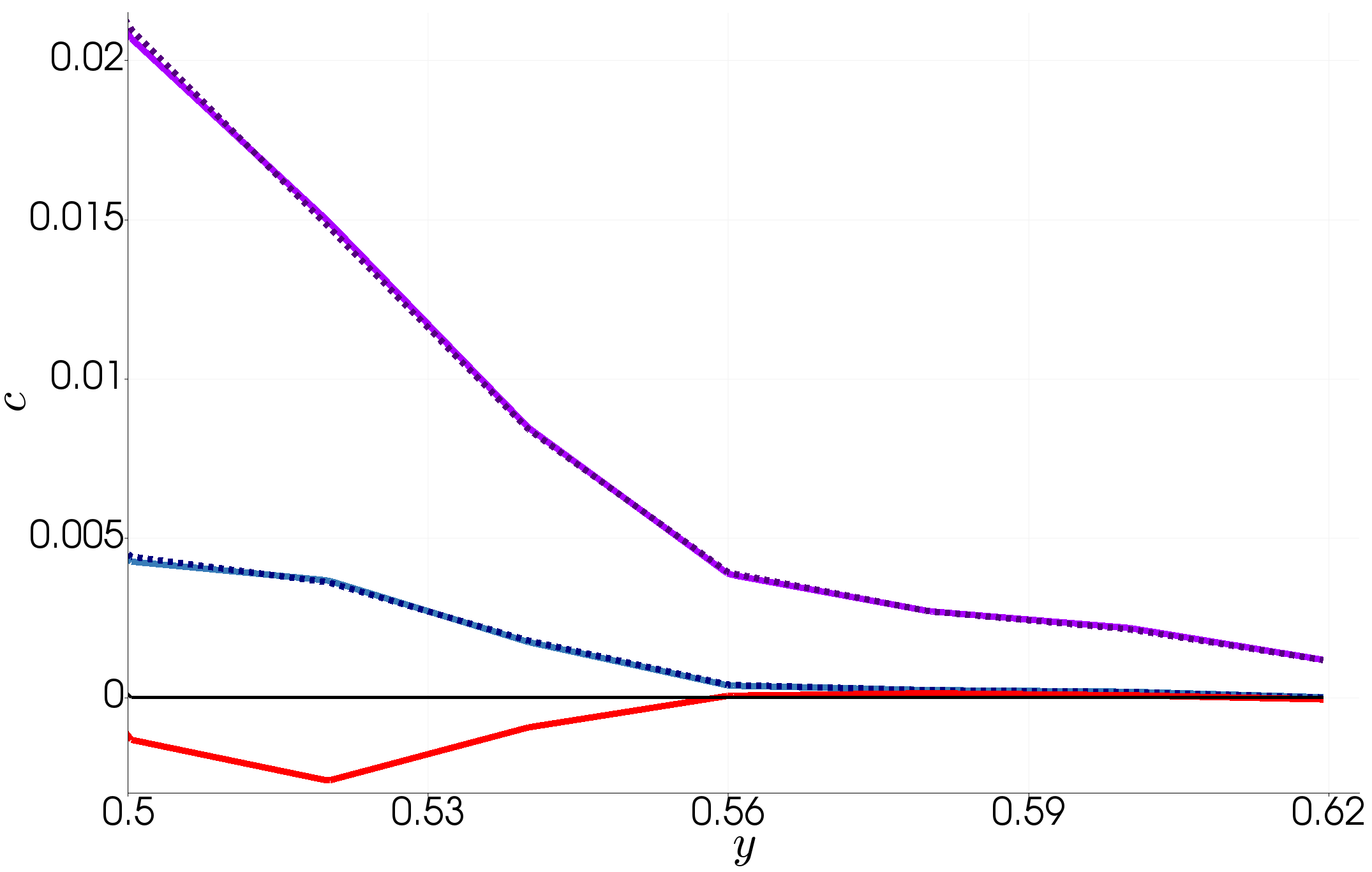}};
\draw[dotted] (4.32,0.6) rectangle +(4.25,0.5);
\draw[dotted] (4.32,0.6) -- (3.5,2.05);
\draw[dotted] (8.57,0.6) -- (8.6,2.1);

\end{tikzpicture}
\captionof{figure}{Comparison of the concentration values for the original model and the transformed model with different DC techniques on a line in the upper part of the outflow ($x=2$, $y>0.4$).}
\label{fig:CompUpperOutflowLine-Channel}
\end{center}
\end{minipage}
\end{figure}

Fig.~\ref{fig:CompWall-Channel} shows a comparison of the concentration on the lower wall of the channel. One can clearly see that the transformation of \eqref{eq:UpperTrans} is very natural and keeps the concentration below 1, while giving results close to the analytical solution.
In order to evaluate the different DC approaches, we compare the different results in Fig.~\ref{fig:CompUpperOutflowLine-Channel} on a line on the upper part of the outflow ($x = 2$, $y > 0.5$), where the analytical solution of the concentration is $c = 0$. The original, untransformed model shows negative concentration values, which can be eliminated by the DC operator. Although the linear form of the discontinuity capturing is able to eliminate the negative values completely, it introduces a much higher amount of numerical diffusion compared to the quadratic form. The results for the isotropic and the crosswind-dissipation DC are almost identical. For the following test cases, we will therefore investigate the behavior of the transformed equation \eqref{eq:AR-transformed} with the DC-quad, CWD method.

\subsection{Two-dimensional Pump}

\begin{figure}
\begin{minipage}{.48\textwidth}
\center\includegraphics[width=\textwidth]{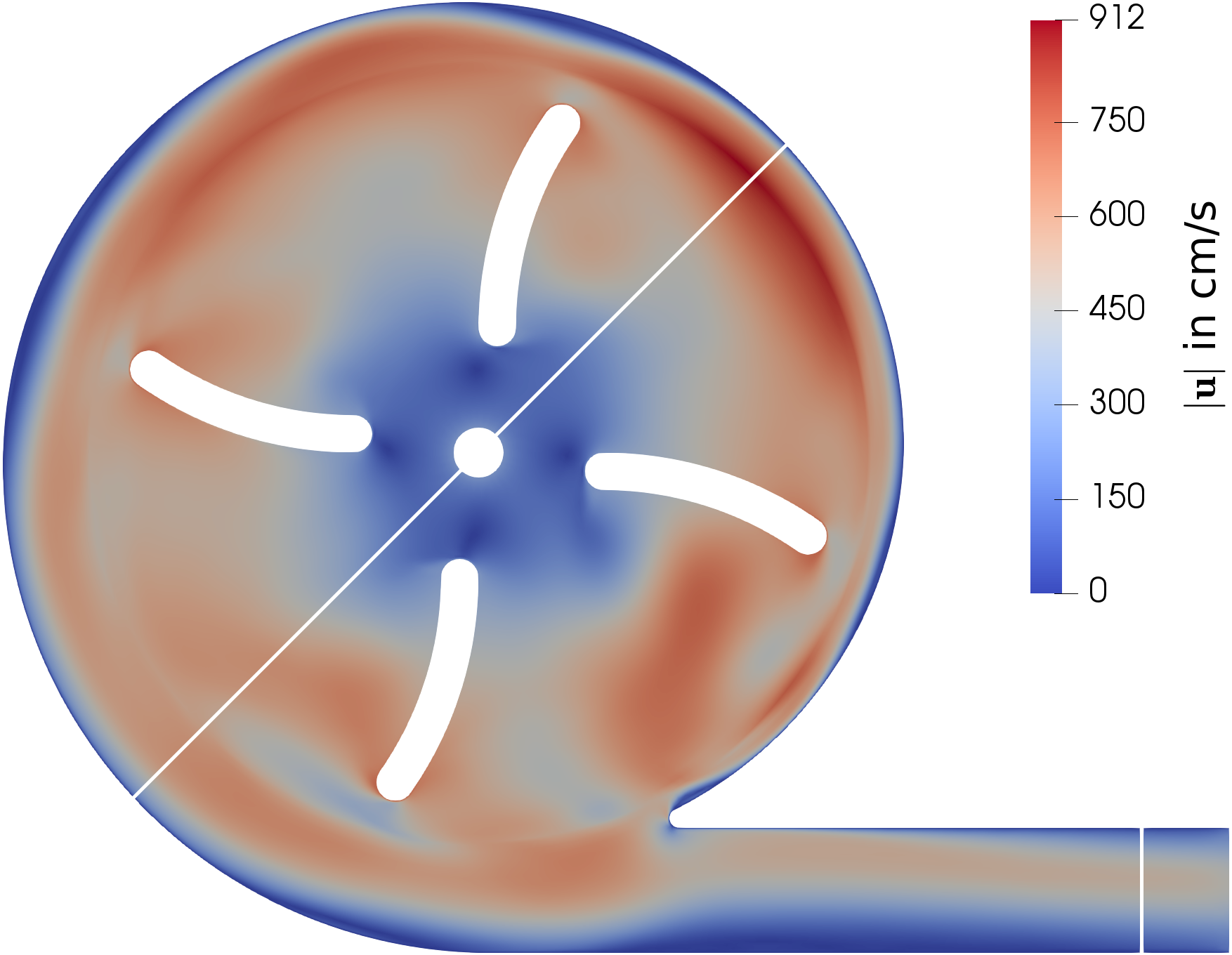}
\caption{Velocity in the 2D pump for $\SI{2500}{rpm}$ computed with the MRF method.\\\null\\\null\\\null}
\label{fig:VelMRF-Pump2D}
\end{minipage}\hfill
\begin{minipage}{.48\textwidth}
\center\includegraphics[width=\textwidth]{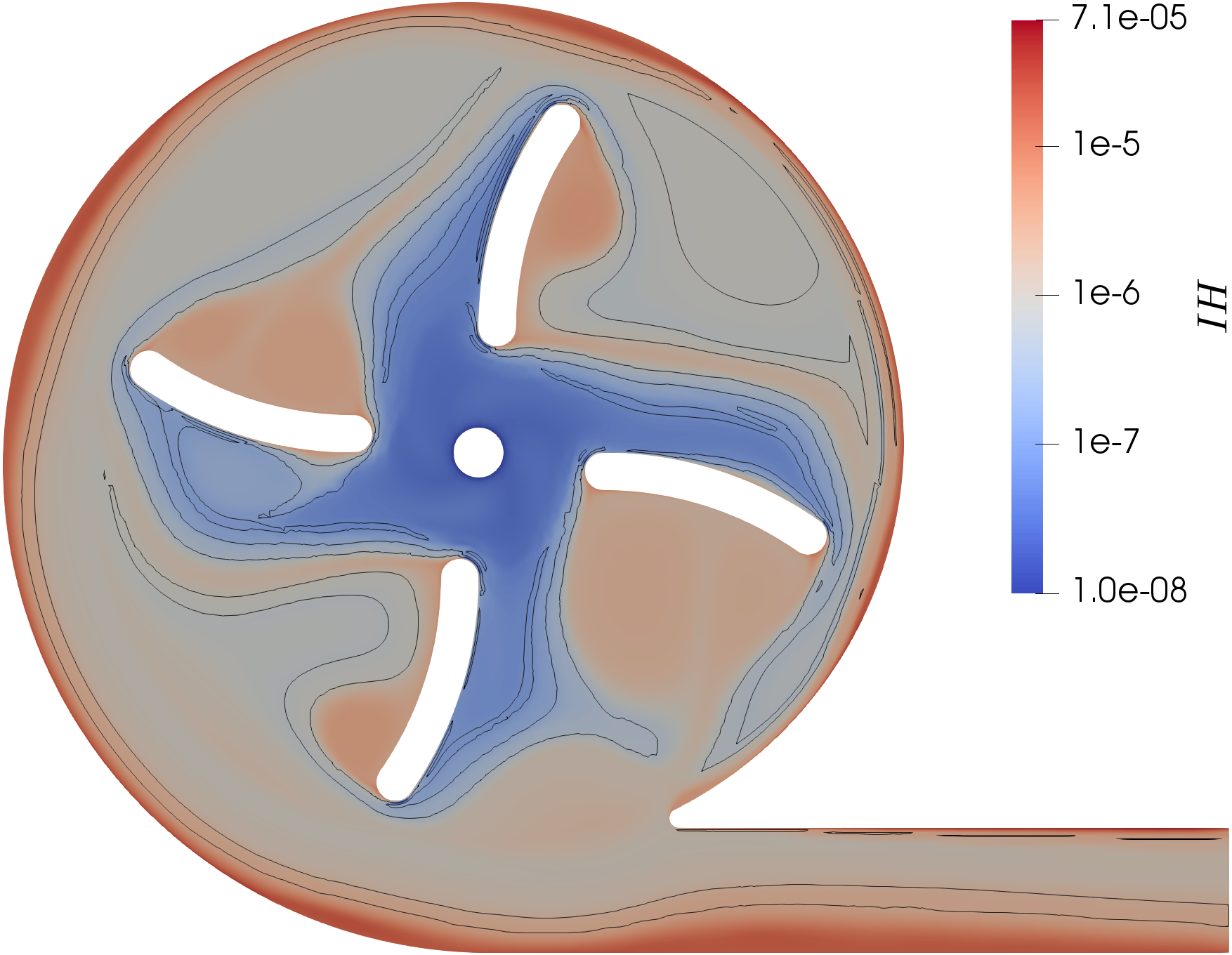}
\captionof{figure}{Hemolysis prediction with the power law model using the transformation~\eqref{eq:UpperTrans} and the DC-quad, CWD method. Additionally, the contours of the areas with negative values for the original model are shown.}
\label{fig:IH-TransDCquadCWD-Pump2D}
\end{minipage}
\end{figure}

\enlargethispage{2mm}As a next test case we consider a two-dimensional pump with a circular inflow in the center of the impeller region and four rounded blades. The computational mesh consists of $n_e = \num{37528}$ triangular elements and $n_n = \num{19287}$ nodes. We computed a quasi-steady blood flow solution using a Newtonian blood model with a density of $\rho = \SI{1.054}{g/cm^3}$ and viscosity of $\nu = \SI{0.035}{g/cm/s}$, utilizing the multiple reference frames (MRF) method. The angular velocity of the impeller is $\SI{2500}{rpm}$ and the inflow velocity is $\SI{500}{cm/s}$ in radial direction. The flow field is shown in Fig.~\ref{fig:VelMRF-Pump2D}.

We compute the stress-based hemolysis for the 2D pump test case using the power law model with the Zhang \etal~\cite{Zhang2011a} parameters: $A = \num{1.228e-7}$, $\alpha = \num{1.9918}$, and $\beta = \num{0.6606}$. We set the linearized hemolysis to 0 at the inflow, assuming unhemolyzed blood is entering the domain. We further make use of the MRF method also for these simulations, to compute a steady concentration field. Since we solve for the linearized hemolysis concentration $l_{IH}$, negative values pose a problem as they would result in complex hemolysis predictions for values of $0 < \beta < 1$. Therefore, we try to minimize negative values by using the crosswind-dissipation DC-quad formulation and set remaining negative values to zero before the conversion to $IH$. The minimal value of $IH = \num{-7.9e-9}$ for the original model can be increased by this approach to $IH = \num{-1.6e-10}$. The areas with negative concentration values for the original model, depicted by their contours in Fig.~\ref{fig:IH-TransDCquadCWD-Pump2D}, are almost completely removed with the DC. Although the concentration values for this test case are far away from the saturation $IH = 1$, we nevertheless make use of the change of variable \eqref{eq:UpperTrans}, due to its simplicity. The resulting hemolysis estimation is shown in Fig.~\ref{fig:IH-TransDCquadCWD-Pump2D}.

\begin{figure}
\begin{minipage}{.6\textwidth}
\begin{tikzpicture}
\node[inner sep=0pt,anchor=south west] at (0,0) {\includegraphics[width=.96\textwidth]{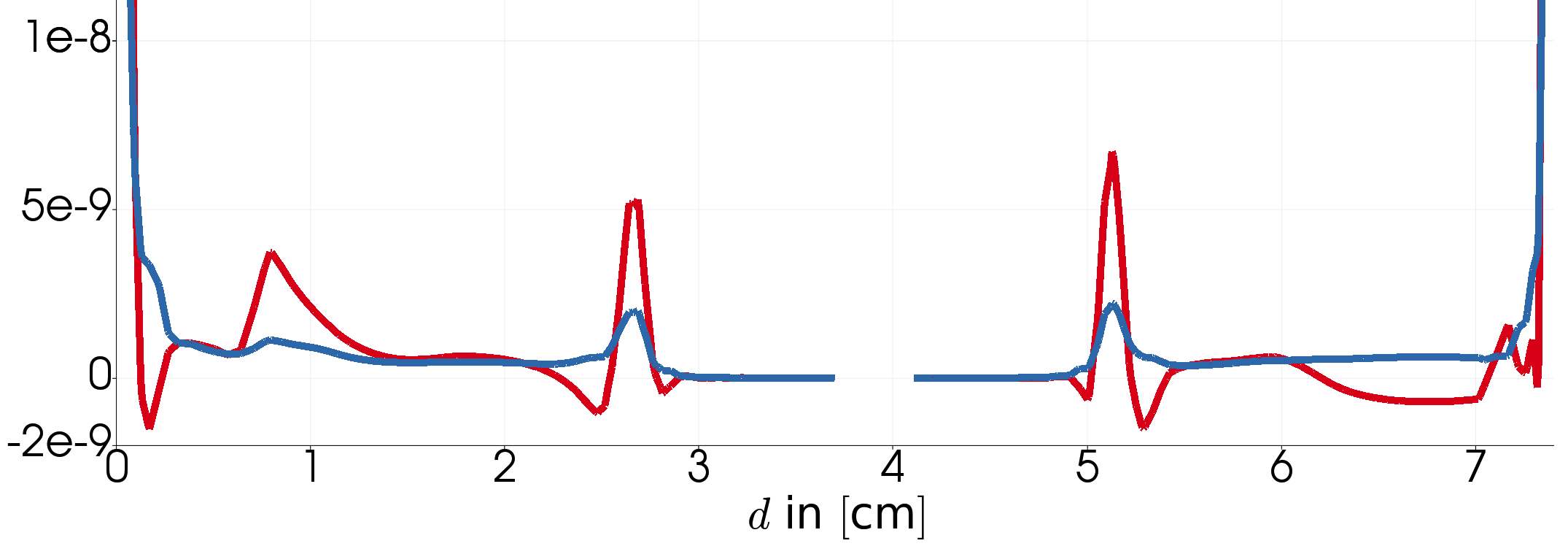}};
\node[inner sep=0pt,anchor=south west] at (0,4) {\includegraphics[width=.96\textwidth]{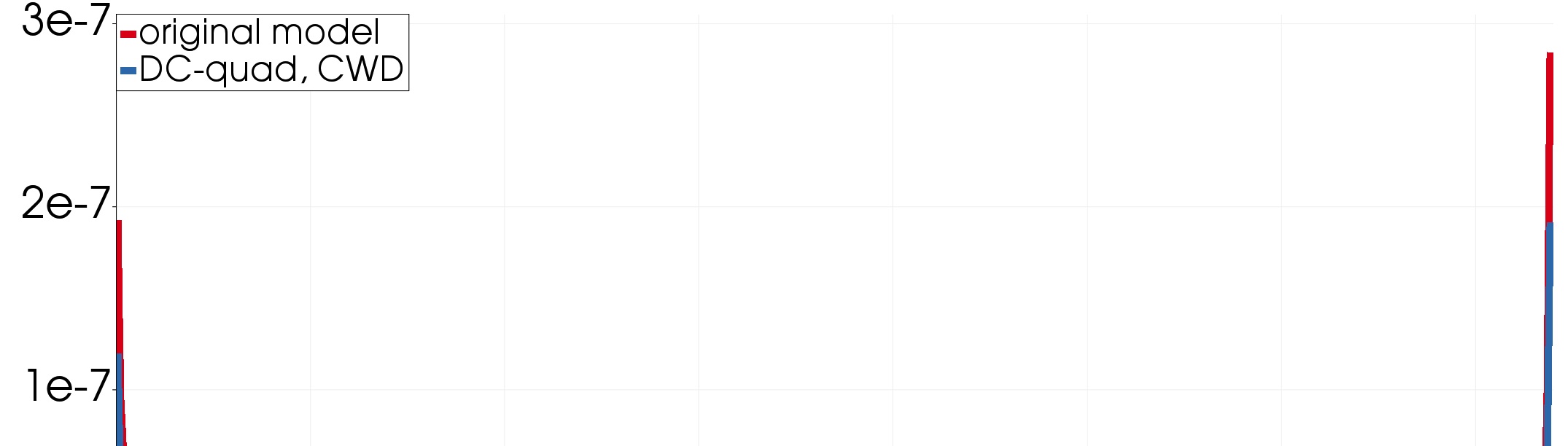}};
\draw (0.8,3.68) -- (0.7,3.77) -- (0.9,3.91) -- (0.8,4);
\node[rotate=90] at (-0.2,3.8) {$l_{IH}$};
\end{tikzpicture}
\captionof{figure}{Comparison of the original model to the transformed model with the DC-quad, CWD method on a diagonal line through the pump chamber (see Fig.~\ref{fig:VelMRF-Pump2D}).}
\label{fig:CompConc-Chamber-Pump2D}
\end{minipage}\hfill
\begin{minipage}{.35\textwidth}
\begin{tikzpicture}
\node[inner sep=0pt,anchor=south west] at (0,0) {\includegraphics[width=.92\textwidth]{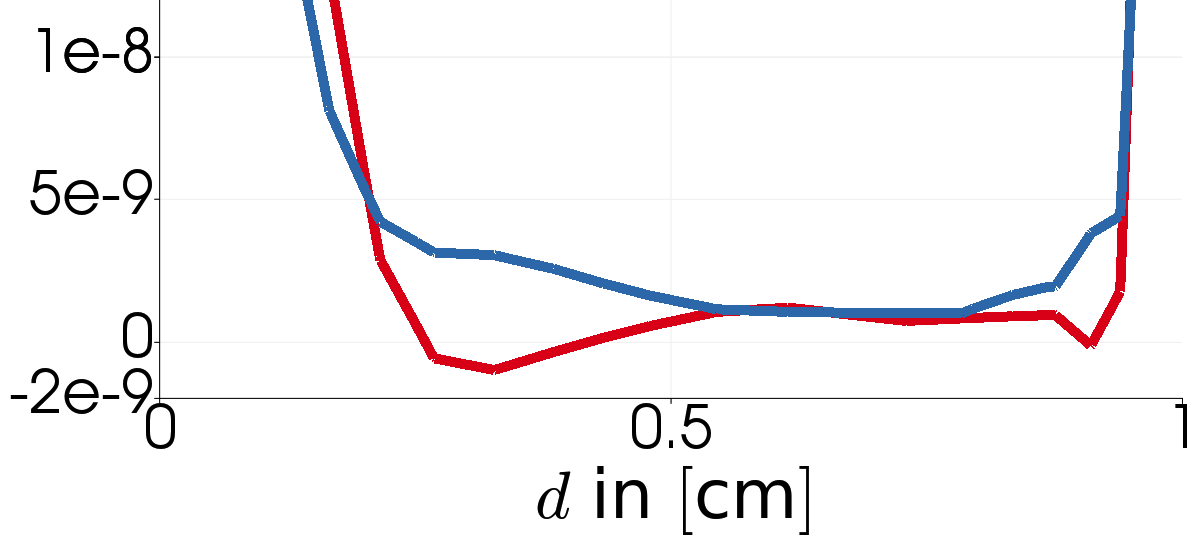}};
\node[inner sep=0pt,anchor=south west] at (0,3) {\includegraphics[width=.92\textwidth]{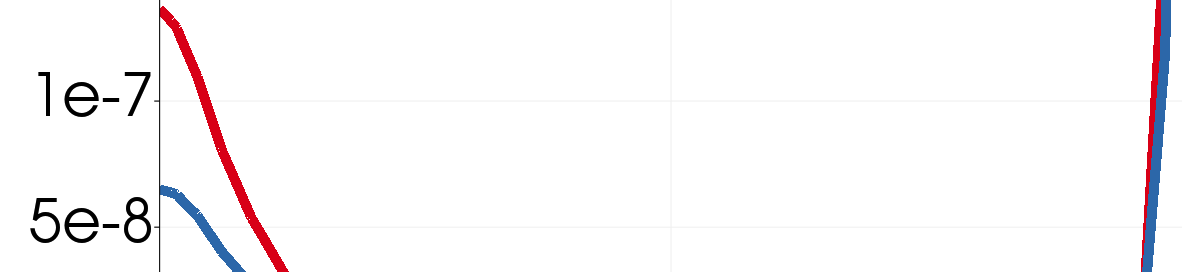}};
\node[inner sep=0pt,anchor=south west] at (0,4.63) {\includegraphics[width=.92\textwidth]{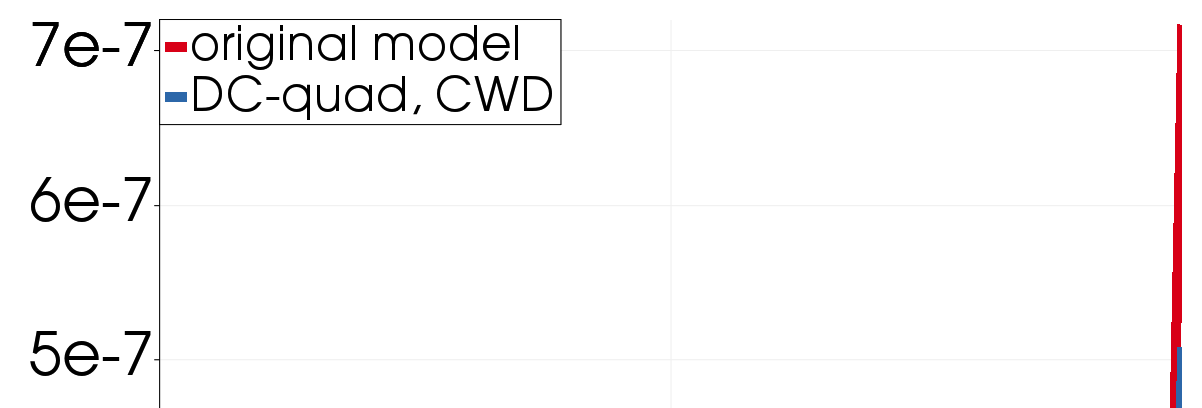}};
\draw (0.79,2.7) -- (0.69,2.78) -- (0.89,2.92) -- (0.79,3);
\draw (0.79,4.33) -- (0.69,4.41) -- (0.89,4.55) -- (0.79,4.63);
\draw[white] (0,7) circle (0.02);
\node[rotate=90] at (-0.2,4) {$l_{IH}$};
\end{tikzpicture}
\captionof{figure}{Same comparison on a line in the outflow tube ($x = \SI{5.3}{cm}$, see Fig.~\ref{fig:VelMRF-Pump2D}).}
\label{fig:CompConc-Outflow-Pump2D}
\end{minipage}
\end{figure}

A comparison of the different estimations is shown in Fig.~\ref{fig:CompConc-Chamber-Pump2D} for a diagonal line through the pump chamber and in Fig.~\ref{fig:CompConc-Outflow-Pump2D} for a line in the outflow tube ($x = \SI{5.3}{cm}$).
It can be seen that the concentration is at least an order of magnitude higher near the pump walls compared to the bulk flow area. Negative $l_{IH}$ values present in the original model predictions are completely eliminated with the DC operator on the lines under consideration. However, the local peak values in the bulk and on the walls are also reduced by the discontinuity capturing (from $IH = \num{8.8e-7}$ to $IH = \num{7.1e-7}$ for the maximum value). In contrast, the investigation of the flow-averaged hemolysis values at the outflow shows a higher concentration for the estimation using the DC. To summarize, the minimization of physically problematic negative values by the discontinuity capturing is very favorable, even if one has to accept alterations in the concentration field, while the correct distribution for this academic test case is unknown, in any case.

\subsection{FDA Benchmark Pump}

The FDA proposed a simple benchmark pump to compare computational flow and hemolysis predictions with experimental data. Particle image velocimetry (PIV) measurements for six operating conditions (Table~\ref{tbl:FDApumpOperatingConditions}) were published by Hariharan \etal~\cite{Hariharan2018a}, and Malinauskas \etal~\cite{Malinauskas2017a} conducted the corresponding hemolysis experiments.
\begin{table}[b!]
\begin{center}
\caption{FDA benchmark blood pump operating conditions.}
\begin{tabular}{clcccc}
\toprule 
      &                  & \multicolumn{4}{c}{Inflow rate} \\
      &                  & $\SI{2.5}{L/min}$ & $\SI{4.5}{L/min}$ & $\SI{6}{L/min}$ & $\SI{7}{L/min}$\\
\midrule
rotor & $\SI{2500}{rpm}$ &         C1        &                   &        C4       & \\
speed & $\SI{3500}{rpm}$ &         C2        &         C3        &        C5       &       C6 \\
\bottomrule
\end{tabular}
\label{tbl:FDApumpOperatingConditions}
\end{center}
\end{table}

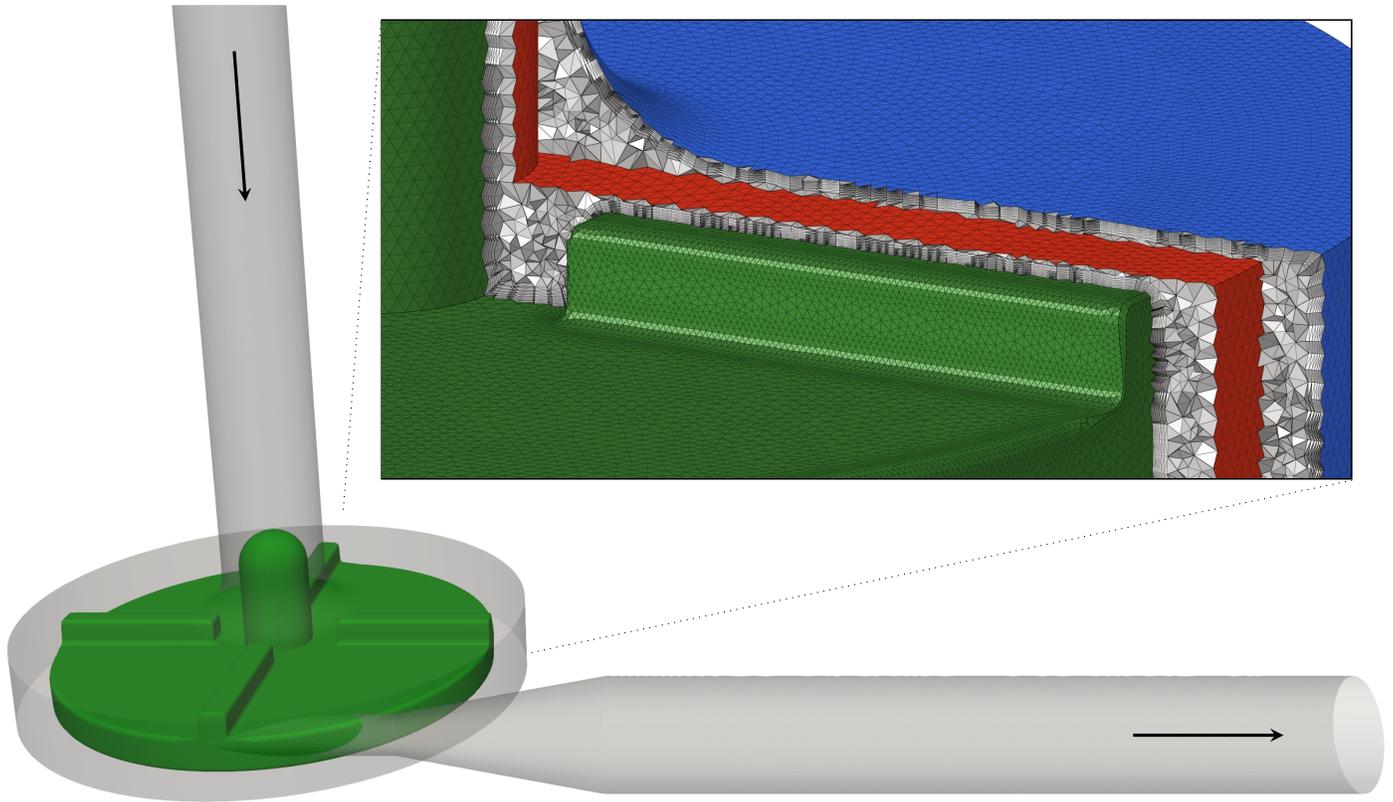
\begin{figure}[p]
\begin{tikzpicture}
\node[inner sep=0pt,anchor=south west] at (0,0) {\includegraphics[width=\textwidth]{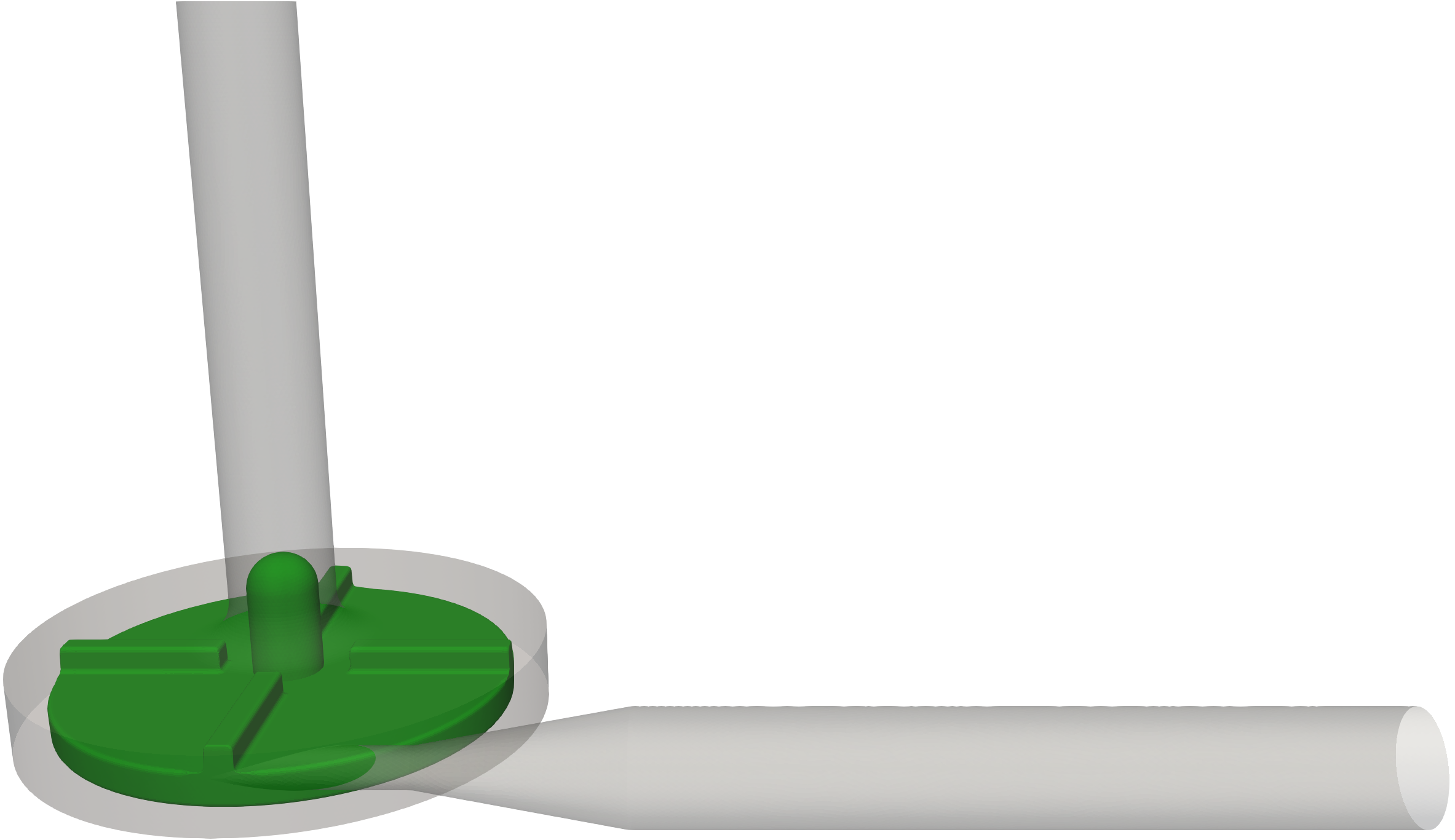}};
\node[inner sep=0pt,anchor=south west,draw,very thick,rectangle] at (5,4.3) {\includegraphics[width=.7\textwidth]{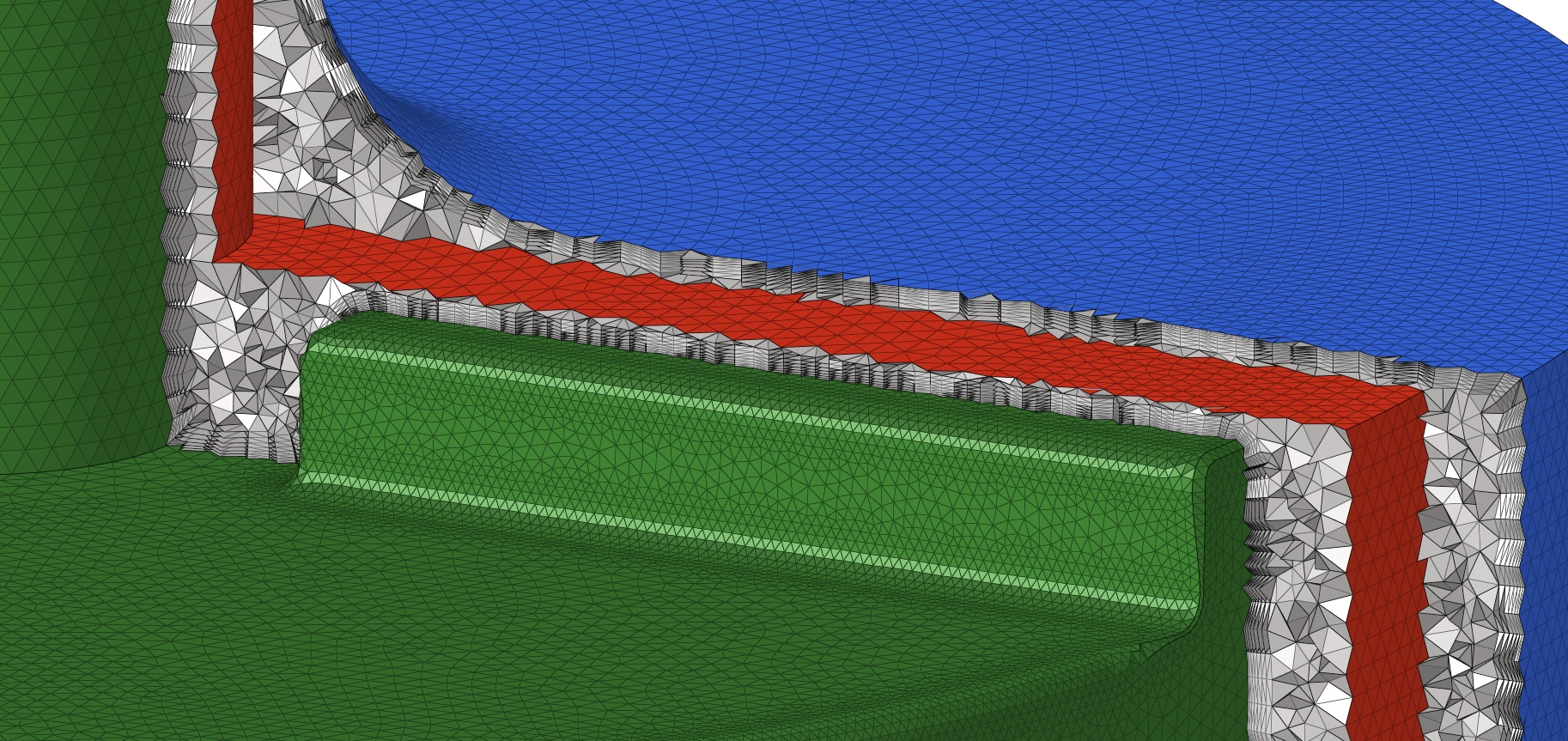}};
\draw[dotted] (7,2) -- (17.95,4.3);
\draw[dotted] (4.5,3.9) -- (5,10.4);
\draw[very thick,-stealth] (3.05,10) -- +(0.15,-2);
\draw[very thick,-stealth] (15,0.9) -- +(2,0);
\end{tikzpicture}
\caption{Geometry and part of the computational mesh for the FDA benchmark blood pump.}
\label{fig:FDApump-mesh}
\end{figure}

In our analysis, we choose a slightly modified geometry of the FDA pump. In order to save computing time, we shortened the inflow and outflow tubes of the pump, since we noticed in previous simulations that their influence on the hemolysis estimation is negligible. We cut the inflow tube $\SI{8}{cm}$ above the pump bottom housing wall and the outflow tube at $z = \SI{12}{cm}$. Furthermore, we did not mesh the rotor shaft, in order to be able to introduce an MRF-interface that is embedded completely in the fluid domain. The computational mesh consists of $n_e = \SI{8.67}{M}$ unstructured tetrahedral elements and $n_n = \SI{1.51}{M}$ nodes and is shown in Fig.~\ref{fig:FDApump-mesh} together with the pump geometry.
\enlargethispage{6mm}In order to capture the velocity gradients near the no-slip walls, we introduce a boundary layer mesh of a total thickness of \SI{500}{\micro m} and five layers with a growth rate of \num{1.2}.

\begin{figure}[p]
\begin{centering}
\begin{tikzpicture}
\node[inner sep=0pt,anchor=south west] at (0,-0.5\textwidth) {\includegraphics[width=.5\textwidth]{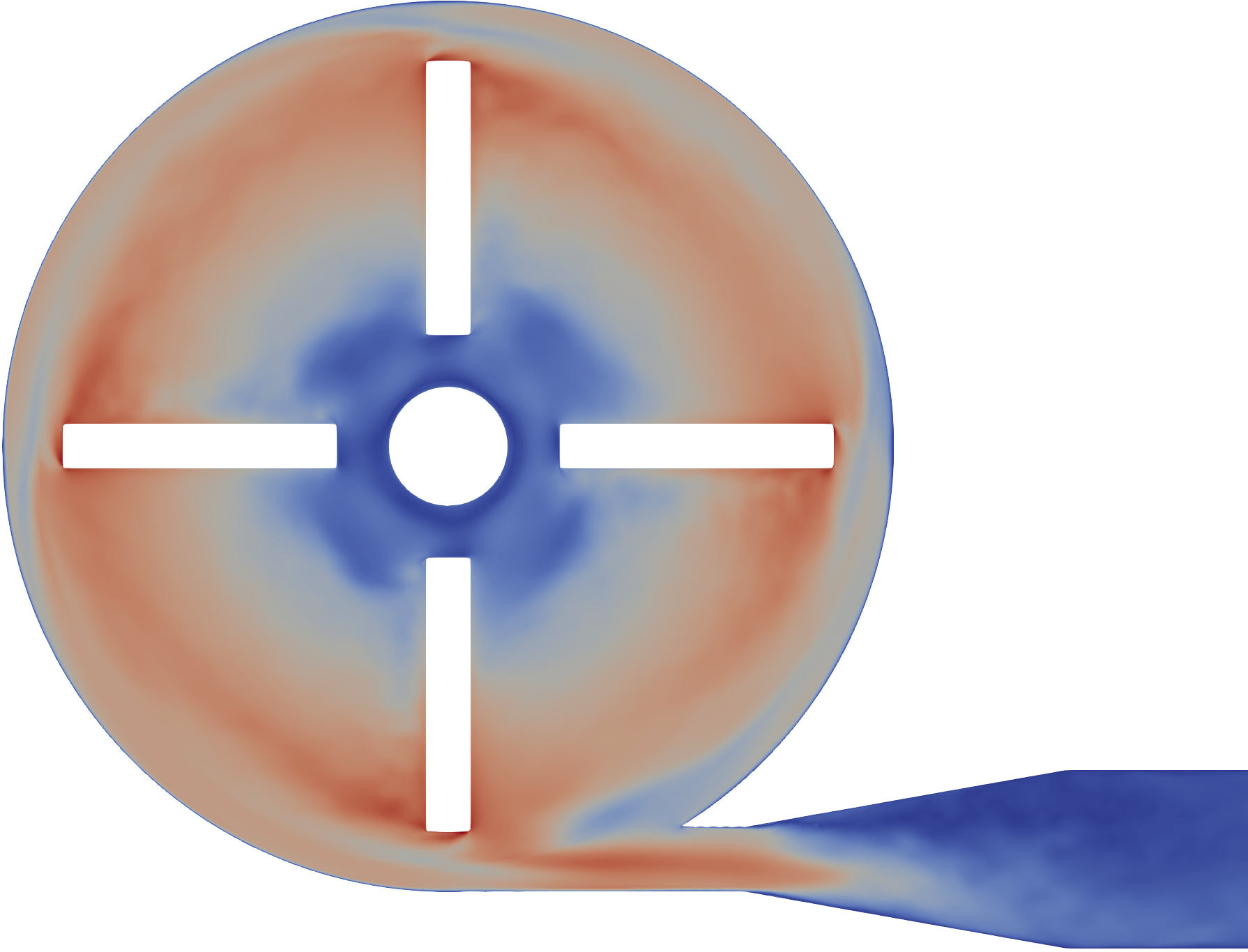}};
\node[inner sep=0pt,anchor=south west] at (.5\textwidth,-0.5\textwidth) {\includegraphics[width=.5\textwidth]{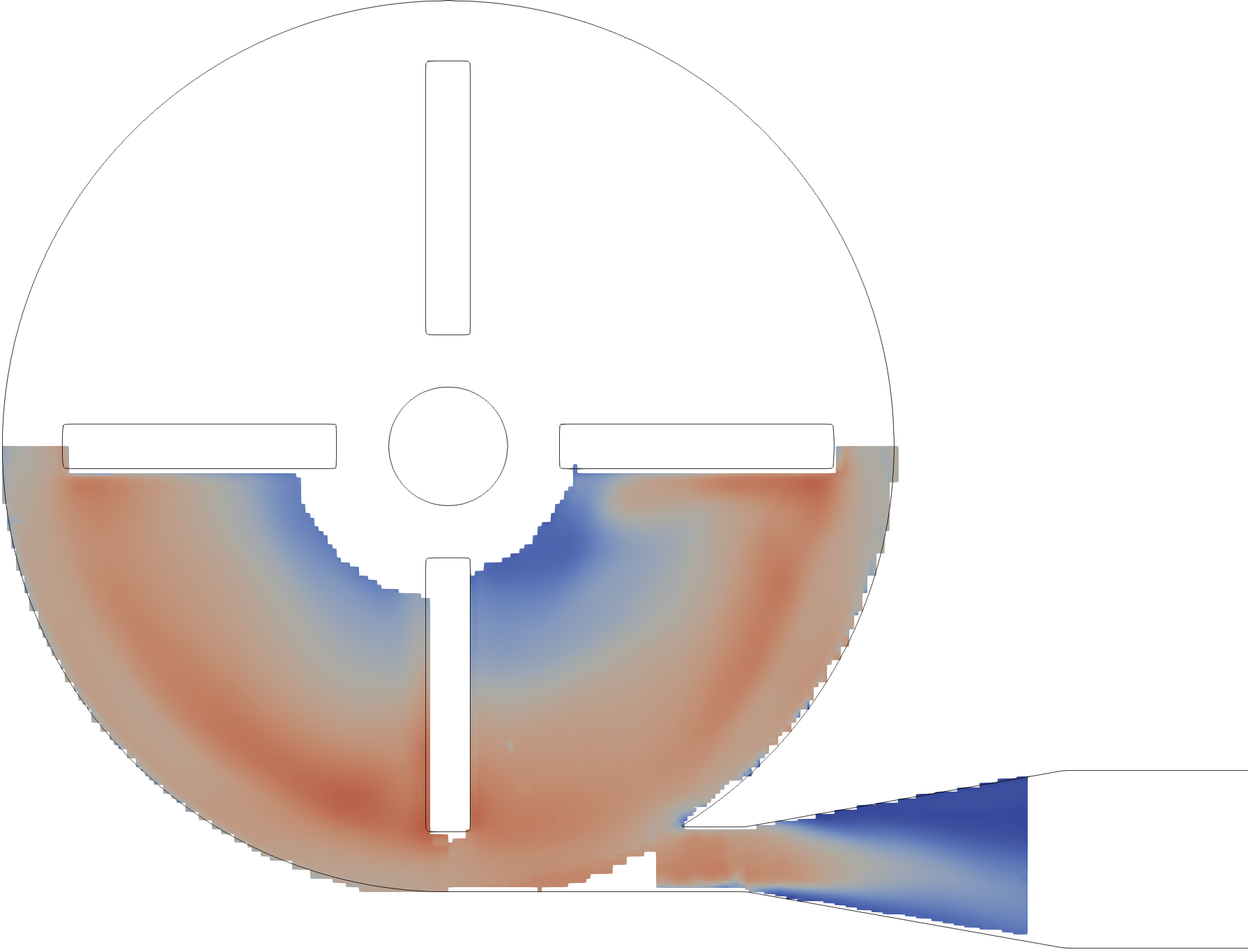}};

\filldraw[color=white] (.5\textwidth,-0.295\textwidth) rectangle +(.48\textwidth,.177\textwidth);

\node[inner sep=0pt,anchor=south west] at (.45\textwidth,-0.2\textwidth) {\includegraphics[width=.4\textwidth]{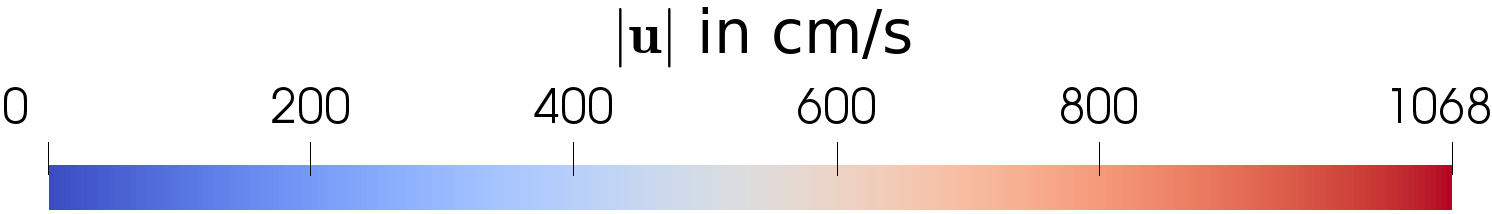}};
\end{tikzpicture}
\end{centering}
\caption{Comparison of CFD results (left) and PIV data \cite{Hariharan2018a} (right) for condition C5.}
\label{fig:FDApump-flow-plane}
\end{figure}
\subsubsection{Flow Simulation}

We compute the flow for the six operating conditions with a Newtonian blood flow model and utilizing the MRF method. For the conditions C1 and C2, we first compute a steady solution for a small Reynolds number, which is used as a restart solution for a transient viscosity ramping computation in $\num{100}$ steps from $\mu = \SI{3.4}{g/cm/s}$ to the desired viscosity of $\mu = \SI{0.034}{g/cm/s}$. From these solutions, we compute transient flow fields for $\num{400}$ time steps and a time step size of $\Delta t = \SI{2.5e-4}{s}$. For the other operating conditions, we restart from the previous condition with identical angular velocity and increase the inflow rate over $\num{40}$ time steps to the desired inflow rate. Again, transient flow fields for $\num{400}$ time steps with a time step size of $\Delta t = \SI{2.5e-4}{s}$ are computed. Finally, these transient solutions are averaged over every $\nth{20}$ time step for the comparison with the PIV measurements \cite{Hariharan2018a} that is presented for condition C5 in Fig.~\ref{fig:FDApump-flow-plane} for an evaluation plane located $\SI{1.2}{mm}$ below the top of the impeller blades.

\begin{figure}
\begin{minipage}{.48\textwidth}
\includegraphics[width=\textwidth]{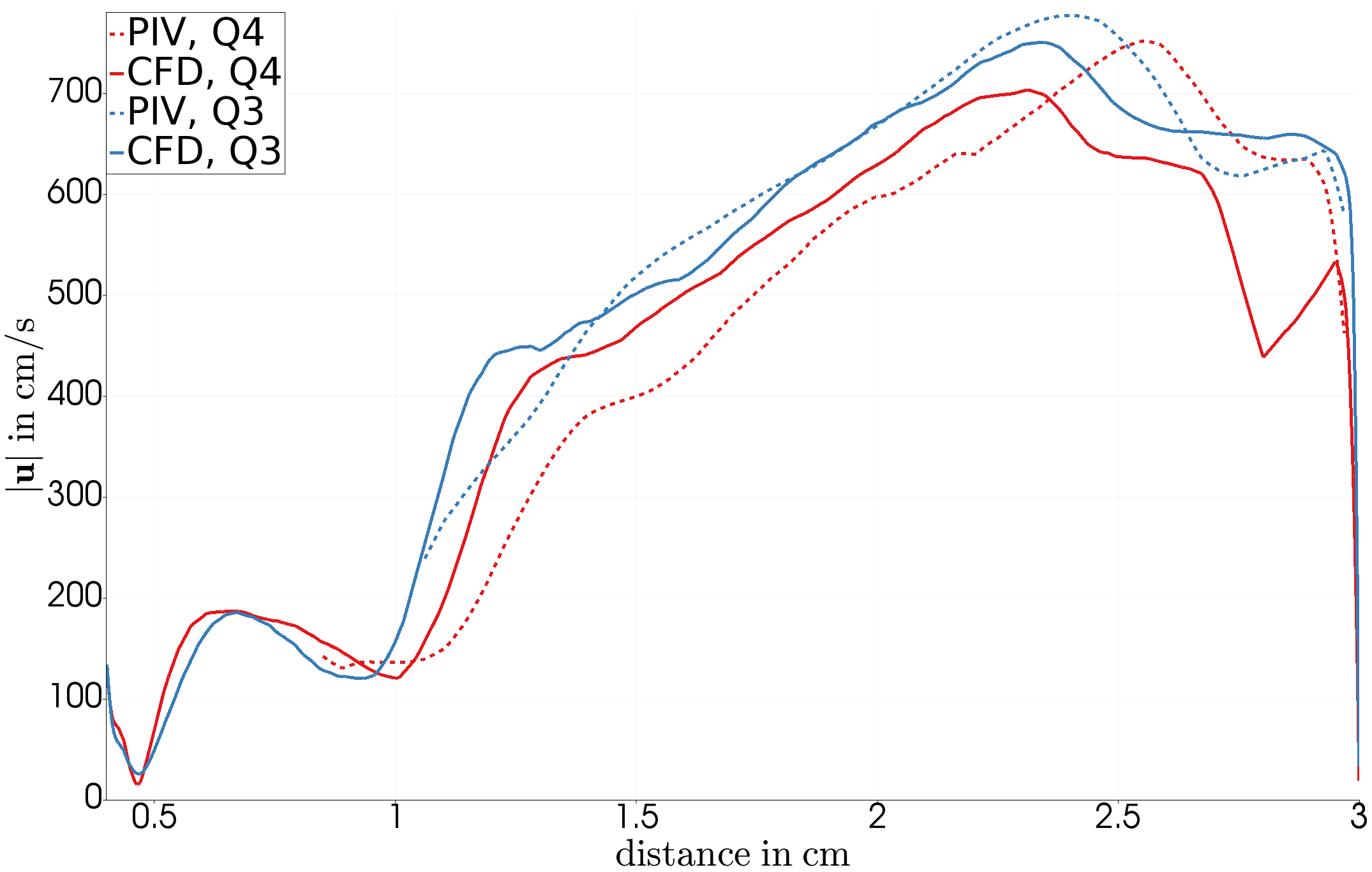}
\captionof{figure}{Comparison of CFD and PIV on two diagonal lines in quadrants Q3 and Q4 in the pump chamber.}
\label{fig:FDApump-flow-line-chamber}
\end{minipage}\hfill
\begin{minipage}{.48\textwidth}
\includegraphics[width=\textwidth]{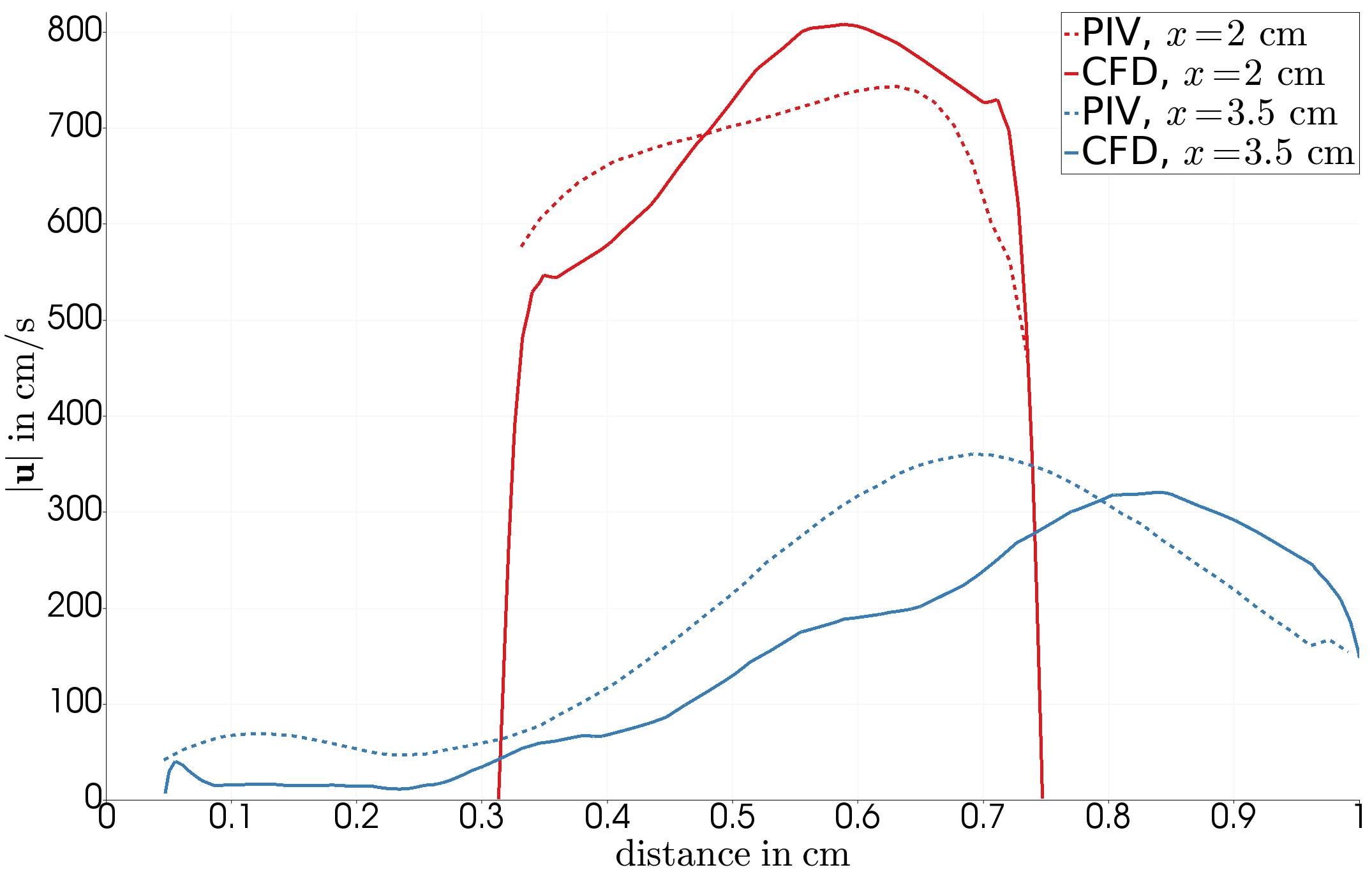}
\captionof{figure}{Comparison of CFD and PIV on two lines in the diffusor region.\\\null}
\label{fig:FDApump-flow-line-diffusor}
\end{minipage}
\end{figure}

Figs.~\ref{fig:FDApump-flow-line-chamber}~and~\ref{fig:FDApump-flow-line-diffusor} show the comparison of the CFD results and the PIV data on different lines in the pump chamber and the diffusor region.
Overall, we find a satisfactory agreement between the simulation results and the experimental data.

\subsubsection{Pore Model Hemolysis Estimation}

\begin{figure}
\begin{minipage}{.48\textwidth}
\includegraphics[width=\textwidth]{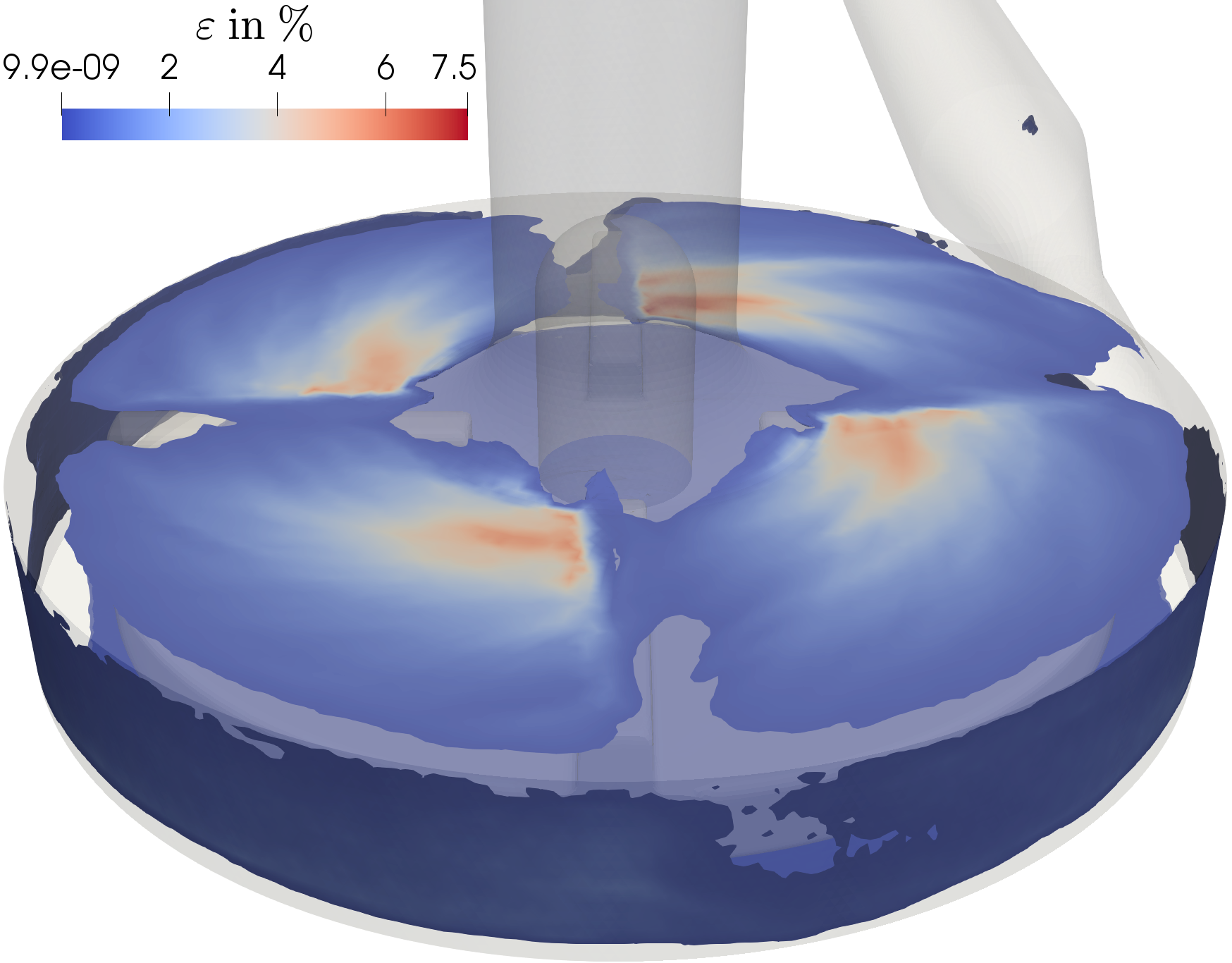}
\captionof{figure}{Surface area strain above the threshold value for condition C5.}
\label{fig:FDA-areaStrain}
\end{minipage}\hfill
\begin{minipage}{.48\textwidth}
\includegraphics[width=\textwidth]{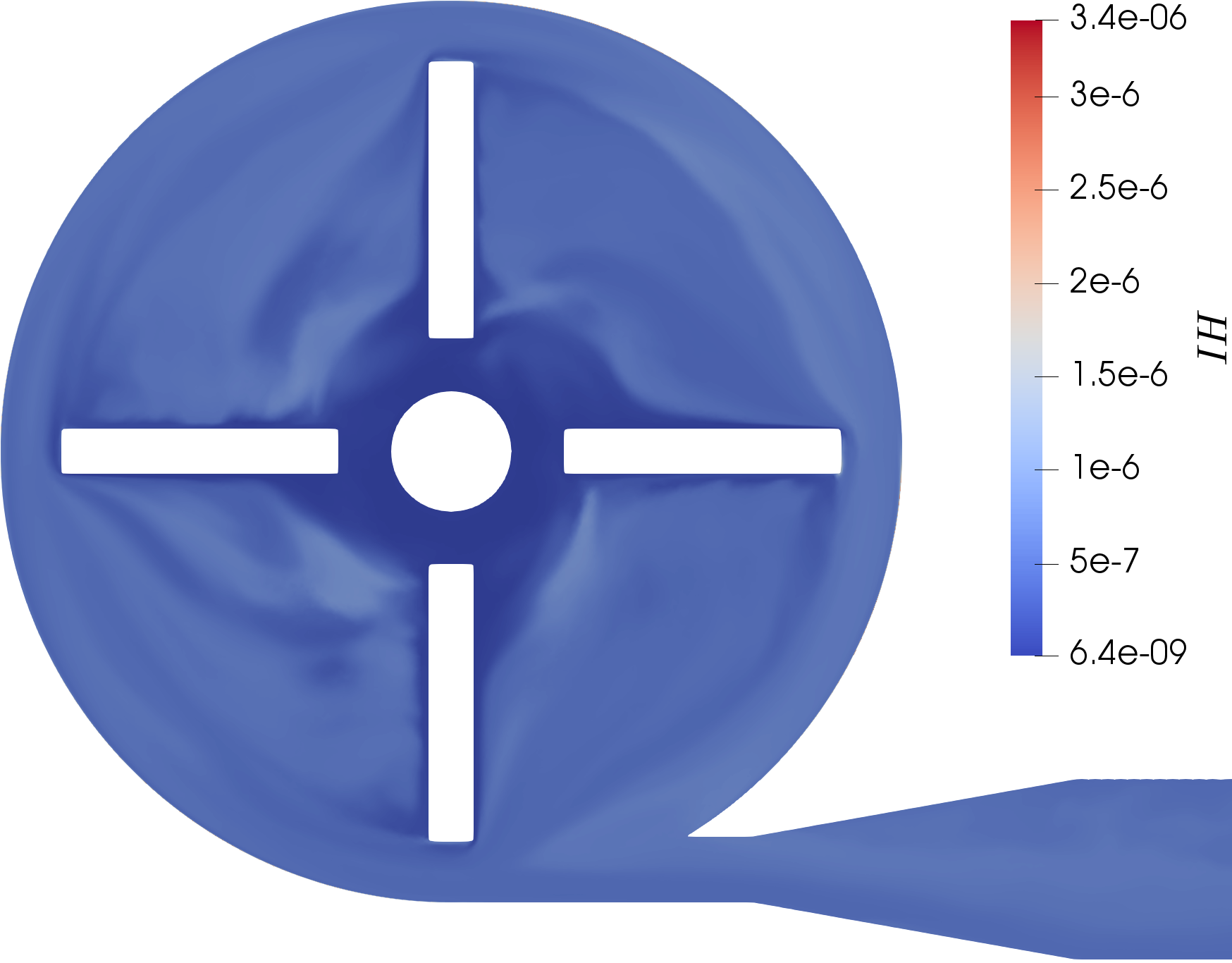}
\captionof{figure}{Index of hemolysis in the evaluation plane for condition C5.}
\label{fig:FDA-IH-plane}
\end{minipage}
\end{figure}

We will use the hemolysis experimental data of Malinauskas \etal~\cite{Malinauskas2017a} as a validation for the hemolysis prediction with the pore model. To this end, we first have to compute the surface area strain of the RBCs for the six flow conditions, utilizing the variational multiscale log-morphology approach presented by us in \cite{Hassler2019a}. We find that the hemolyzing areas, with surface area strains above the threshold value, are mainly situated close to and at the pump chamber walls, which is shown in Fig.~\ref{fig:FDA-areaStrain} for condition C5.

The next step is to compute the generation and distribution of free plasma hemoglobin with the pore model \eqref{eq:poreModel}. The corresponding flow fields and surface area strain distributions for the six operating conditions are the pre-computed inputs for the model. At the inflow, we assume that unhemolyzed blood is entering the domain by setting the concentration to $\bar{c}_\mathrm{in} = 0$. We use the parameters $(h,k) = (\num{4.48e-8},\num{1.31})$ for porcine blood and utilize the concentration transformation \eqref{eq:UpperTrans} and the crosswind-dissipation discontinuity capturing operator with the quadratic diffusion definition. Furthermore, we use the MRF method to compute a quasi-steady hemoglobin distribution. \enlargethispage{1mm}This distribution for condition C5 in the evaluation plane is shown in Fig.~\ref{fig:FDA-IH-plane}.

\begin{figure}
\begin{center}
\begin{tikzpicture}
\begin{axis}[
    width = .6\textwidth,
    height = .4\textwidth,
    ybar,
    legend style={at={(0.17,0.94)},
      anchor=north,legend columns=1},
    xlabel={operating condition},
    ylabel={$\Delta PHb$ in $\si{mg/dL}$},
    ytick={0,50,100,150,200,250,300},
    symbolic x coords={C1,C2,C3,C4,C5,C6},
    xtick=data,
    bar width=15pt,
    ]
\addplot+ [line width=1pt] coordinates {(C1,16.82) (C2,41.71) (C3,39.43) (C4,13.19) (C5,36.11) (C6,33.97)};
\addplot+ [
    line width=1pt,
    error bars/.cd,
    y dir = plus,y explicit,
    error bar style={line width=2pt},
    error mark options={rotate=90,line width=2pt,mark size=4pt}]
    coordinates {
        (C1,21) +- (0,28)
        (C2,29) +- (0,37)
        (C3,53) +- (0,39)
        (C4,45.5) +- (0,30.5)
        (C5,44) +- (0,25)
        (C6,156) +- (0,118)};
\legend{pore model,experiment}
\end{axis}
\end{tikzpicture}
\belowcaptionskip = -5mm
\caption{Comparison of the estimated and measured free plasma hemoglobin.}
\label{fig:FDA-poreModel-Validation}
\end{center}
\end{figure}
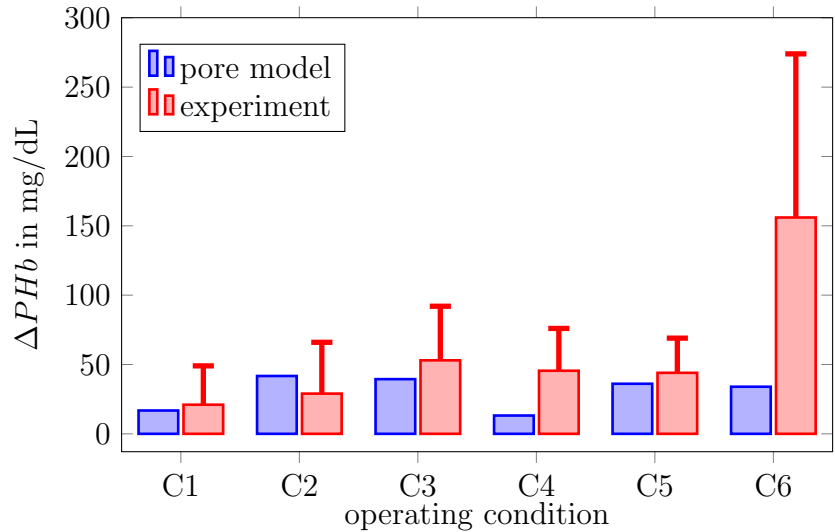

In order to compare the computational results with the experimental findings, we have to convert the flow-averaged hemolysis index at the outflow $\Gamma_\mathrm{out}$ (surface normal $\boldsymbol{n}$)
\begin{equation}
\overline{IH}_\mathrm{out} = \frac{\int_{\Gamma_\mathrm{out}} \bu\cdot\boldsymbol{n} \, IH \, \diff \Gamma}{\int_{\Gamma_\mathrm{out}} \bu\cdot\boldsymbol{n} \, \diff \Gamma}
\end{equation}
to a free plasma hemoglobin content $\Delta PHb$. By our choice of the inflow boundary condition, we assume that $\overline{IH}_\mathrm{out}$ is proportional to the modified index of hemolysis, and hence, the conversion can be done by
\begin{equation}
\Delta PHb = \overline{IH}_\mathrm{out}\frac{Hb}{1-Hct}\frac{Q\cdot T}{V_\mathrm{loop}},
\end{equation}
with the hemoglobin content of blood $Hb$, the blood hematocrit value $Hct = \SI{36}{\%}$, the flow rate $Q$, the experiment duration $T = \SI{120}{min}$, and the volume of the blood sample in the mock loop $V_\mathrm{loop} = \SI{250}{mL}$ for the Malinauskas \etal~\cite{Malinauskas2017a} experiments. Since they do not report a hemoglobin content, we assume a value of $Hb = \SI{15000}{mg/dL}$. The validation of our simulation predictions shown in Fig.~\ref{fig:FDA-poreModel-Validation} shows a very satisfactory agreement with the experimental results.
Similar computations with the stress-based and strain-based model, relying on the power law, give hemolysis results that are about two orders of magnitude higher than the experimental values. This emphasizes the superior behavior of the physiologically motivated pore model for hemolysis estimations.

\section{Conclusion}

We presented a change of variable approach and discontinuity capturing techniques to restrict the numerical simulation of advection-reaction equations to physical ranges. To our best knowledge, our proposed change of variable \eqref{eq:UpperTrans} is new for the advection-reaction equation. It is able to fulfill the constraint of an upper bound on the concentration field without introducing another level of complexity to the equation. Due to its general form and its simplicity, it can be applied to other discretization techniques and might be useful for similar models (e.g., advection-diffusion-reaction eqs.).

The discontinuity capturing operator that we proposed to use is a combination of the DC operator of Shakib~\etal\cite{Shakib91b}, defined on the reference element, and the crosswind-dissipation approach of Codina~\cite{Codina93a}. We also used the numerical diffusion definitions of Shakib \etal~\cite{Shakib91b} and linearized them by using the concentration values from the last time step, in order to achieve convergence of the simulations. Although this method also introduces alterations to the positive concentration values compared to a solution without a DC, it proves very useful due to its elimination of nonphysical negative values for the most part.

We showed the performance of these methods for the hemolysis prediction for two academic test cases. For the FDA benchmark blood pump, we successfully applied the proposed methods for a pore formation hemolysis model and found a good agreement with the experimental values.

\subsection*{Acknowledgments}

\enlargethispage{6mm}The authors gratefully acknowledge support from the German Research Foundation (DFG) grant BE~3689/15 "Drug-Eluting Coronary Stents in Stenosed Arteries: Medical Investigation and Computational Modelling". We especially thank Dr. Lutz Pauli for fruitful discussions and guidance, and Maria Torrentallé Dot for her numerous contributions to the pore model during her master thesis; both of them made this article possible. The authors gratefully acknowledge the computing time granted through JARA-HPC on the supercomputer JURECA at Forschungszentrum Jülich.

\bibliography{hemo.bib}
\bibliographystyle{fine}

\end{document}